\documentclass{article}
\usepackage[utf8]{inputenc}
\usepackage{amsmath}
\usepackage{amsfonts}
\usepackage{amssymb}
\usepackage{amsthm}
\usepackage{soul}
\usepackage{bm}
\usepackage[shortlabels]{enumitem}
\usepackage{comment}
\usepackage{graphicx}
\usepackage{float}
\usepackage{hyperref}
\usepackage{arydshln}
\usepackage{lipsum}

\newtheorem*{Theorem 1}{Theorem 1.1}
\newtheorem{Lemma}{Lemma}
\newtheorem*{Lemma 6.1}{Lemma 6.1}

\newtheorem{Theorem}{Theorem}
\newtheorem{Proposition}{Proposition}
\newtheorem{Corollary}{Corollary}

\newcommand{\norm}[1]{\left\lVert#1\right\rVert}
\newcommand{\Longupdownarrow}{\Big\Updownarrow}
\newcommand{\rectangle}{\fboxsep0pt\fbox{\rule{1.3em}{0.5pt}\rule{1pt}{1.2ex}}}

\newcommand\blfootnote[1]{%
  \begingroup
  \renewcommand\thefootnote{}\footnote{#1}%
  \addtocounter{footnote}{-1}%
  \endgroup
}

\newsavebox{\overlongequation}
\newenvironment{dontbotheriftheequationisoverlong}
 {\begin{displaymath}\begin{lrbox}{\overlongequation}$\displaystyle}
 {$\end{lrbox}\makebox[0pt]{\usebox{\overlongequation}}\end{displaymath}}

\title{\textbf{Constructing Displacement Vectors}}
\author{Alon Agin}
\date{}

\begin{document}
\maketitle
\begin{abstract}
 \ \\ Let $\vec{v}\in\mathbb{R}^2\setminus\mathbb{Q}^2$ ,let $\norm{\cdot}$ be an arbitrary norm on $\mathbb{R}^2$, and let $\big(q_n,\vec p_n\big)_{n=0}^{\infty}
\\ \subset \mathbb{N} \times \mathbb{Z}^{2}$ be the best approximation vectors sequence of $\vec{v}$ with respect \\ to $\norm{\cdot}$. We define the nth long displacement vector of $\vec{v}$ to be \\ 
$\vec{\beta}_n:=\sqrt{q_{n+1}} \, (q_{n}\vec{v}-\vec{p}_n)$ and prove the existence of long displacement vectors who have non-typical properties – focusing on their length, direction, and congruence class.   
\end{abstract}

\blfootnote{The author wishes to thank Barak Weiss for his invaluable guidance, support, and comments on earlier versions of the paper.
The author also thanks Nikolay Moshchevitin and Renat Akhunzhanov for helpful discussions. This paper represents the author’s M.Sc. dissertation at Tel Aviv University under the supervision of Barak Weiss. Financial support from the grants BSF 2016256, ISF 2019/19 and ISF-NSFC 3739/21 is gratefully acknowledged.}

\section{Introduction}
\ \\ Throughout this paper bold symbols and texts are used for definitions. \\ For $\vec{v} \in \mathbb{R}^{d} \setminus \mathbb{Q}^{d}$ and arbitrary norm $\lVert \cdot \lVert$,
denote $\boldsymbol{\langle\vec{v} \rangle_{\lVert \cdot \lVert}}$ := min $ \{ \, \lVert \vec{v}-\vec{p} \lVert \,\, \mid \vec{p}\in\mathbb{Z}^{d} \, \}$. 

\noindent Using Minkowski's Convex Body Theorem, one deduces a generalised Dirichlet Theorem – there exists a minimal positive constant $\gamma_{\lVert \cdot \lVert}$ such that for all $1<T\in\mathbb{R}$ there exists $q\in\mathbb{N}$ with $1\leq q\leq T$ such that \,$T^{\frac{1}{d}} \, \langle q \vec{v} \rangle_{\lVert \cdot \lVert}
\leq\gamma_{\lVert \cdot \lVert}$. 

\noindent One then easily concludes that the inequality  \,$q^{\frac{1}{d}} \, \langle q \vec{v} \rangle_{\lVert \cdot \lVert}
\leq\gamma_{\lVert \cdot \lVert}$ has infinitely many solutions for $q\in\mathbb{N}$. 

\ \\ Define \textbf{the best approximation vectors sequence of} \boldsymbol{$\vec{v}$} (with respect to $\lVert \cdot \lVert$) to  be  a  sequence
$ \big(q_n,\vec p_n\big)_{n=0}^{\infty}=\big(q_n(\vec v),\vec p_n(\vec v)\big)_{n=0}^{\infty}
\subset \mathbb{N} \times \mathbb{Z}^{d} $ \,\,such that: 
\begin{itemize}
\item $1=q_0<q_1<...<q_n<...$
\item $q_n=\min\,\,\{q\in\mathbb{N}\,\mid\,\langle q\vec{v}\rangle_{\lVert \cdot \lVert}<\langle q_{n-1}\vec{v}\rangle_{\lVert \cdot \lVert}\}$
\item $\vec{p}_n$ realises $\langle q_n\vec{v}\rangle_{\lVert \cdot \lVert}$.  I.e.  $\langle q_n\vec{v}\rangle_{\lVert \cdot \lVert} = \lVert q_n\vec{v}-\vec{p}_n \lVert.
$
\end{itemize} 

\newpage \noindent Notice that the vector $\vec{p}_n$ is not uniquely defined. Nevertheless, by discreteness of $\mathbb{Z}^d$, this ambiguity is relevant only for (at most) a finite number of indices. Hence, all best approximations sequences agree starting from some index. 

\noindent Consequently, we refer to all of these sequences as if they were unique.

\ \\ Going back to Dirichlet Theorem, we now conclude that for all $n \geq 0$ 
\begin{equation}
\label{equation1}
   q_{n}^{\frac{1}{d}}\,\, \lVert q_n\vec{v}-\vec{p}_n \lVert  \, < q_{n+1}^{\frac{1}{d}}\,\,   \lVert q_n\vec{v}-\vec{p}_n \lVert  \, \leq \, \gamma_{\lVert \cdot \lVert}. 
\end{equation} \ \\ This leads us to define natural objects for the study of Diophantine approximations. If $\big(q_n,\vec p_n\big)_{n=0}^{\infty}$ is the best approximations sequence of $\vec{v}$ (w.r.t $\norm{\cdot}$) then define

\vspace*{0.15cm}
\textbf{the nth short displacement vector} := $q_{n}^{\frac{1}{d}} \, (q_{n}\vec{v}-\vec{p}_n)$ 
  
\vspace*{0.15cm}  
\textbf{the nth long displacement vector} := $q_{n+1}^{\frac{1}{d}} \, (q_{n}\vec{v}-\vec{p}_n)$.

\ \\ Note that the long and short displacement vectors have the same direction, but different length. 

\ \\  Many of the classical Diophantine approximations results obtained in the last century can be phrased in terms of these displacement vectors. We denote 
$$
\boldsymbol{\vec{\alpha}_n(\vec{v})}:=\vec{\alpha}_n(\vec{v},\norm{\cdot})=q_n^{\frac{1}{d}} \, (q_{n}\vec{v}-\vec{p}_n)
$$
$$\boldsymbol{\vec{\beta}_n(\vec{v})}:=\vec{\beta}_n(\vec{v},\norm{\cdot})=q_{n+1}^{\frac{1}{d}} \, (q_{n}\vec{v}-\vec{p}_n)$$ \ \\ keeping in mind these vectors are norm dependent.

\ \\ For the case $d=1$, properties of the sequences of the short and long displacement vectors have been studied by many, especially in the context of continued fractions. For higher dimensions, Davenport and Schmidt showed in \cite{DS} that for the maximum norm $\norm{\cdot}=\norm{\cdot}_{\infty}$, for Lebesgue almost every $\vec{v} \in \mathbb{R}^{d}$ we have $\limsup\limits_{n\rightarrow\infty} \,\lVert\,\vec{\beta}_n(\vec{v})\lVert_\infty=1$. This result has been generalised for arbitrary norms gradually, culminating with Kleinbock and Rao who recently showed in \cite{KR1} that for an arbitrary norm $\norm{\cdot}$, for almost every $\vec{v} \in \mathbb{R}^{d}$ we have  $\limsup\limits_{n\rightarrow\infty} \,\lVert\,\vec{\beta}_n(\vec{v})\lVert=\gamma_{\norm{\cdot}}$ ("almost all vectors in $\mathbb{R}^d$ are not Dirichlet-improvable"). For the lower bound, Chevallier showed in \cite{Chev1} that for almost every $\vec{v}\in\mathbb{R}^d$ we have $\liminf\limits_{n\rightarrow\infty} \,\lVert\,\vec{\beta}_n(\vec{v})\lVert=0$.

 \ \\ Regarding direction, Rogers observed in \cite{Rog} that for the maximum norm, two consecutive displacement vectors cannot lie in the same quadrant.  
 
 \noindent Moshchevitin showed in \cite{Mosh} that this result cannot be generalised to arbitrary norms. 

\ \\ For the case of $d=1$, Bosma, Jager and Wiedijk showed in \cite{BJW} that there exist probability measures $\mu_1,\mu_2$ on $\mathbb{R}$ such that for Lebesgue almost every $v\in\mathbb{R}$ the sequence of short displacement vectors $(\alpha_n(v))_{n=0}^{\infty}$ equidistributes with respect to $\mu_1$, and the sequence of long displacement vectors $(\beta_n(v))_{n=0}^{\infty}$ equidistributes with respect to $\mu_2$. This result has been generalised recently by Shapira and Weiss in \cite{SW} for an arbitrary dimension and norm for the case of the short displacement vectors (with the measure $\mu_1$ being norm dependent). Also in \cite{SW}, Shapira and Weiss generalised a distribution theorem by Moeckel \cite{Moe} regarding congruence properties of the sequence $(q_n,\vec{p}_n)$. Informally, they showed that for a fixed $m\in\mathbb{N}$, the $(d+1)$ integer tuples  $\,(q_n,p_{n,1},...\,,p_{n,d})$ are distributed uniformly between residues modulo $m$ which are "valid" (for example need to have $\gcd\,(q_n,p_{n,1},...\,,p_{n,d})=1$), independently of the norm under consideration. 

\ \\ So we have some understanding of how a typical sequence of the form $(\vec{\beta}_n(\vec{v}))_{n=0}^{\infty}$ behaves. What about non-typical behaviour?  \\ Define \textbf{the Dirichlet spectrum of $\mathbb{R}^d$} (with respect to $\norm{\cdot}$) to be
\begin{align*}
\boldsymbol{\mathbb{D}_{d,\norm{\cdot}}}:=&\, \{\,\limsup\limits_{T\rightarrow\infty} \,T^{\frac{1}{d}} \,\bigl( \min_{1\leq q\leq T} \langle q \vec{v} \rangle_{\lVert \cdot \lVert}\bigr)
\mid \vec{v}\in\mathbb{R}^{d} \setminus\mathbb{Q}^{d}\,\} \ \\ 
=&\,\{\,\limsup\limits_{n\rightarrow\infty} \,\lVert\,\vec{\beta}_n(\vec{v})\lVert\,\,\mid\,\vec{v}\in\mathbb{R}^{d} \setminus\mathbb{Q}^{d}\,\}.
\end{align*} By (\ref{equation1}) we know $\mathbb{D}_{d,\norm{\cdot}}\subseteq[0,\gamma_{\norm{\cdot}}]$. For the case $d=1$, by many works on continued fractions, we know that the spectrum has quite a complicated structure which we still don't understand completely (see e.g \cite{AS}). \\ In \cite{AS}, Akhunzhanov and Shatskov were the first to consider higher dimensions. Surprisingly, they showed that for $d=2$ and for the Euclidean norm $\norm{\cdot}=\norm{\cdot}_2$ we have $\mathbb{D}_{2,\norm{\cdot}_2}=[0,\gamma_{\norm{\cdot}_2}]=[0,\sqrt{\frac{2}{\sqrt{3}}}]$. Recently, Schleischitz showed in \cite{Sch} that for any $2\leq d$  and for the maximum norm $\norm{\cdot}=\norm{\cdot}_\infty$ we have $\mathbb{D}_{d,\norm{\cdot}_\infty}=[0,\gamma_{\norm{\cdot}_\infty}]=[0,1]$. He also showed that if $\norm{\cdot}$ is an expanding norm on $\mathbb{R}^d$ then there exists a constant $M_{\norm{\cdot}}$ such that $[0,M_{\norm{\cdot}}]\subseteq \mathbb{D}_{d,\norm{\cdot}} $, where we say \textbf{a norm is expanding} if for all $\vec{s}\in\mathbb{R}^d$ and all $1\leq j\leq d$ we have that $\norm{\vec{s}}\geq\norm{\pi_j(\vec{s})}$, where $\pi_j$ is the orthogonal projection to the $j$th coordinate axis. Also recently, Kleinbock and Rao showed in \cite{KR2} that for $d=2$ and an arbitrary norm, $\gamma_{\norm{\cdot}}$ is an accumulation point of $\mathbb{D}_{2,\norm{\cdot}}$. 

\ \\ In \cite{AM}, Akhunzhanov and Moshchevitin introduced a more delicate notion than the Dirichlet spectrum. We say that $\vec{v}$ is \textbf{badly approximable} if $\inf\limits_{T\geq1} \,T^{\frac{1}{d}} \,\bigl( \min_{1\leq q\leq T} \langle q \vec{v} \rangle_{\lVert \cdot \lVert}\bigr)>0$. \ \\ Now define \textbf{the BA Dirichlet spectrum of} $\boldsymbol{\mathbb{R}^d}$ (with respect to $\norm{\cdot}$) to be
\begin{align*}
\boldsymbol{\mathbb{D}_{d,\norm{\cdot}}^{BA}}:=&\, \{\,\limsup\limits_{T\rightarrow\infty} \,T^{\frac{1}{d}} \,\bigl( \min_{1\leq q\leq T} \langle q \vec{v} \rangle_{\lVert \cdot \lVert}\bigr)
\mid \vec{v}\,\,is\,\,badly\,\, approximable\,\,\} \ \\ 
=&\,\{\,\limsup\limits_{n\rightarrow\infty} \,
\lVert\,\vec{\beta}_n(\vec{v})\lVert\,\,\mid \vec{v}\,\,is\,\,badly\,\, approximable\,\,\} \ \\ 
=&\,\{\,\limsup\limits_{n\rightarrow\infty} \,
\lVert\,\vec{\beta}_n(\vec{v})\lVert\,\,\mid \,\, \inf\limits_{n\geq1}\,\lVert\vec{\alpha}_n(\vec{v})\lVert\,>0 
\,\}.
\end{align*}\ \\ Clearly we have $\mathbb{D}_{d,\norm{\cdot}}^{BA}\subseteq\mathbb{D}_{d,\norm{\cdot}}$. In \cite{AM} it was shown that for $d=2$ and the Euclidean norm $\lVert\cdot\lVert=\lVert\cdot\lVert_2$ we have that $\mathbb{D}^{BA}_{2,\norm{\cdot}}$ is dense in $[0,\sqrt{\frac{2}{\sqrt{3}}}]$.

\ \\ \ \\ We suggest another refinement of the Dirichlet spectrum. We say that \textbf{$\boldsymbol{\vec{v}}$ is Dirichlet converging} (w.r.t) $\norm{\cdot}$) if $\lim\limits_{n\rightarrow\infty} \,\lVert\vec{\beta}_n(\vec{v})\lVert$ exists, and define \textbf{the \\limiting Dirichlet spectrum of $\boldsymbol{\mathbb{R}^d}$} (with respect to $\norm{\cdot}$) to be
$$
\boldsymbol{\mathbb{D}_{d,\norm{\cdot}}^{\lim}}:=\,\{\,\lim\limits_{n\rightarrow\infty} \,
\,\lVert\vec{\beta}_n(\vec{v})\lVert
\,\,\mid \, \vec{v}\,\,is\,\,Dirichlet\,\,converging\,\}.
$$
\ \\ As far as we know, very little is known regarding Dirichlet converging vectors and the limiting Dirichlet spectrum. Again we have $\mathbb{D}_{d,\norm{\cdot}}^{\lim}\subseteq\mathbb{D}_{d,\norm{\cdot}}$.  

\ \\ In \cite{AS}, Akhunzhanov and Shatskov also proved that  $\mathbb{D}_{2,\norm{\cdot}_2}^{\lim}=[0,\sqrt{\frac{2}{\sqrt{3}}}]$. As far as we know, this is the only result concerning $\mathbb{D}_{d,\norm{\cdot}}^{\lim}$ for $d\geq2$. \,Furthermore, as far as we know, \cite{AS} is the first and only proof that the set of Dirichlet converging vectors is even an infinite set (for the case $d=2, \norm{\cdot}=\norm{\cdot}_2$). 

\ \\ In \cite{AS}, Akhunzhanov and Shatskov actually proved a stronger statement than just $\mathbb{D}_{2,\norm{\cdot}_2}^{\lim}=\mathbb{D}_{2,\norm{\cdot}_2}=[0,\sqrt{\frac{2}{\sqrt{3}}}]$. They showed that for the case \\ $d=2, \norm{\cdot}=\norm{\cdot}_2$ one can have a much finer control on the lengths of the long displacement vectors rather than only its limit -- controlling it at each step. Explicitly, they proved that if $(I_n)_{n=0}^{\infty}$ is a sequence of closed, non-degenerate intervals with $I_n\subset[0,\sqrt{\frac{2}{\sqrt{3}}}\,]$, then there exists continuum many $\vec{v}\in\mathbb{R}^2$ such that $\lVert\,\vec{\beta}_n(\vec{v})\lVert_2\,\in I_n$ for all $n$.

\ \\ In this paper we generalise techniques developed by Akhunzhanov and Shatskov in \cite{AS} in order to show that for an arbitrary norm on $\mathbb{R}^2$, one can construct continuum many vectors $\vec{v}$ such that the sequence  $(\vec{\beta_n}(\vec{v}))_{n=0}^{\infty}$ has non-typical properties regarding length, direction, and congruence class. 

\newpage \noindent For $\theta$,$\delta$ with $0\leq\theta\leq2\pi$ and $0<\delta$ we define
$$\boldsymbol{\Omega(\theta,\delta)}:= \{\,(\,r\cos(\beta),r\sin(\beta)\,) \in\mathbb{R}^{2} \mid r\in\mathbb{R}\,\,and\,\, \theta-\delta\leq\beta\leq\theta+\delta\,\}.$$

 \ \\ We say that \textbf{$\boldsymbol{\norm{\cdot}_@}$ is a best ellipsoid norm for $\boldsymbol{\norm{\cdot}}$ with parameters $\boldsymbol{a,c>0,D\geq1}$}\, if for all $(u_1,u_2)\in\mathbb{R}^2$ it satisfies \begin{equation}\label{equation2}
\norm{\vec{u}}_@:=\sqrt{au_1^2+cu_2^2} 
\, \leq \, \lVert \vec{u} \lVert \, 
\leq \, D\sqrt{au_1^2+cu_2^2}=D\norm{\vec{u}}_@
\end{equation}\ \\and if it maximises the value \begin{equation}\label{equation3}
 M_{\lVert \cdot \lVert}:= \left(\frac{4ac}{4D^2-1}\right)^{\frac{1}{4}}.
\end{equation} 

\ \\  In this paper we prove the following:

\begin{Theorem}
\begin{flushleft}
\ \\
\ \\
Let $\lVert \cdot \lVert$ be an arbitrary norm on $\mathbb{R}^2$. Then there exist an explicit constant $M_{\lVert \cdot \lVert}>0$ such that the following holds:

\ \\ Let $\theta$,$\delta$ with $0\leq\theta<2\pi$ and $\,0<\delta$.

\vspace*{0.4cm}

Let $(I_{n})_{n=0}^{\infty}$ be a sequence of closed non-degenerate intervals with $I_{n}\subset[0,M_{\lVert \cdot \lVert}]$ for all n.

\vspace*{0.4cm}

Let $m\in\mathbb{N}$, and let $(z_{n})_{n=1}^{\infty}\in\{0,1,...,m-1\}$ be a sequence of residues mod $m$ such that $z_n\in\{sz_{n-1}+iz_{n-2}\pm z_{n-3}\mod{m}\,|\,\,i,s\in\mathbb{Z}\}$ for all $n$, where $z_0=1$, $z_{-1}=z_{-2}=0$.
\end{flushleft}
\vspace*{0.4cm}
Then there exist continuum many $\vec{v} \in \mathbb{R}^2$ with best approximations sequence \, $(q_n,\vec p_n)_{n=0}^{\infty}$= 
$(q_n(\vec v),\vec p_n(\vec v))_{n=0}^{\infty}$ such that the following hold:

\begin{enumerate}[(A)]
    
    \item \label{1A} \underline{length:}\,\,\,$\lVert\,\vec{\beta}_n(\vec{v})\lVert\,\, \in I_{n}$ for all $0\leq n$.
   
    \item \label{1B} \underline{direction:}\,\, $\vec{\beta}_n(\vec{v}) \in \Omega(\theta,\delta)$ for all $0\leq n$.\\ Furthermore, two consecutive displacement vectors always lie in opposite quadrants  
    
    \item \label{1C} \underline{arithmetical control:}\, $q_n\equiv z_{n} \, mod\,\, m$ for all $1\leq n$.
    \end{enumerate} 
 
\ \\ Furthermore, $\vec{v}$ is not badly approximable, $\{(q_{n-2},\vec{p}_{n-2}),(q_{n-1},\vec{p}_{n-1}),(q_{n},\vec{p}_{n})\}$ form a basis of the lattice $\mathbb{Z}^3$ for all $n\geq2$, and if $\norm{\cdot}_@$ is a best ellipsoid norm for $\norm{\cdot}$ then $M_{\lVert \cdot \lVert}$ is the same as in (\ref{equation3}).

\newpage   In particular,
 \begin{equation}
 \label{equation4}
 \end{equation}

\begin{center}

\begin{tabular}{|c| |c| |c|} 
 \hline
 norm & best ellipsoid parameters & $M_{\norm{\cdot}}$ \\ [0.5ex] 
 \hline\hline
 Euclidean & $b=0,\,a=c=1,\,D=1$ & $\sqrt{\frac{2}{\sqrt{3}}}$  \\ 
 \hline
 maximum & $b=0,\,a=c=\frac{1}{2},\,D=\sqrt{2}$ & $7^{\frac{-1}{4}}$  \\
 \hline
 p-norm ($2<p$) & $b=0,\,a=c=2^{\frac{2}{p}-1},\,D=2^{\frac{1}{2}-\frac{1}{p}}$ & $8^{\frac{1}{2p}}\,(8-4^{\frac{1}{p}})^{\frac{-1}{4}}$   \\
 \hline
  p-norm ($1\leq p<2$) & $b=0,\,a=c=1,\,D=2^{\frac{1}{p}-\frac{1}{2}}$ & $\sqrt{2}\,\,(2^{1+\frac{2}{p}}-1)^{\frac{-1}{4}}$   \\
  [1ex] 
 \hline
\end{tabular}
\end{center}   
\end{Theorem}

\ \\ \textbf{Remark 1.} By equation (\ref{equation1}) we have $M_{\lVert \cdot \lVert}\leq\gamma_{\norm{\cdot}}$, and so $M_{\lVert \cdot \lVert}$ is best possible if $\lVert \cdot \lVert$ is induced from an inner product induced by a diagonal matrix; when $\norm{(x,y)}=\sqrt{ax^2+cy^2}$. For example, if $\lVert \cdot \lVert$ is the standard Euclidean norm then indeed $M_{\lVert \cdot \lVert}=\gamma_{\lVert \cdot \lVert}=\sqrt{\frac{2}{\sqrt{3}}}$. 

\ \\ In short, this follows from the equality 
$\gamma_{\norm{\cdot}}=2\,\,\sqrt{\delta_{\norm{\cdot}}\mathbb{B}^{-1}}$, where $\delta_{\norm{\cdot}}$ is the packing density of $\mathbb{R}^2$ with respect to $\norm{\cdot}$ and $\mathbb{B}$ is the volume of the unit ball of $\norm{\cdot}$. As the packing density of $\mathbb{R}^2$ with respect to the Euclidean norm is $\frac{\pi}{\sqrt{12}}$, the packing density is invariant under linear mappings, and the area of the ellipse defined by $\{(u_1,u_2)\,|\,au_1^2+cu_2^2=1\}$ is equal to $\frac{2\pi}{\sqrt{4ac}}$, we get that if $\norm{(x,y)}=\sqrt{ax^2+cy^2}$\, then $\gamma_{\norm{\cdot}}=(\frac{4ac}{3})^{\frac{1}{4}}$. Plugging $D=1$ in (\ref{equation3}), we get $M_{\lVert \cdot \lVert}=\gamma_{\norm{\cdot}}$.

\ \\ \textbf{Remark 2.} We phrase result (\ref{1C}) in terms of the sequence of denominators $(q_n)_{n=1}^{\infty}$, but one can formulate a similar statement regarding any one of the two coordinates of the sequence $(\vec{p}_n)_{n=1}^{\infty}$.

\ \\ \textbf{Remark 3.} As both $\{(q_{n-2},\vec{p}_{n-2}),(q_{n-1},\vec{p}_{n-1}),(q_{n},\vec{p}_{n})\}$ and \\$\{(q_{n-3},\vec{p}_{n-3}),(q_{n-2},\vec{p}_{n-2}),(q_{n-1},\vec{p}_{n-1})\}$ form a basis of the lattice $\mathbb{Z}^3$, we have that $(q_n,\vec{p}_n)=s(q_{n-1},\vec{p}_{n-1})+i(q_{n-2},\vec{p}_{n-2})\pm(q_{n-3},\vec{p}_{n-3})$ for some $s,i\in\mathbb{Z}$, hence the condition $z_n\in\{sz_{n-1}+iz_{n-2}\pm z_{n-3}\mod{m}\,|\,\,i,s\in\mathbb{Z}\}$ is optimal within this kind of a construction. It
would be interesting to prove a similar statement with $q_n$ being any residue mod $m$, but to do so one cannot have that every 3 consecutive best approximation vectors form a basis of $\mathbb{Z}^3$.

\begin{Corollary} The Dirichlet and the limiting Dirichlet spectrum of $\mathbb{R}^2$ (with respect to arbitrary norm) always contain a segment. \\ In detail, if  $\norm{\cdot}$ is an arbitrary norm on $\mathbb{R}^2$, then there exist $M_{\norm{\cdot}}>0$ such that
$[0,M_{\norm{\cdot}}]\subseteq\mathbb{D}^{\lim}_{2,\norm{\cdot}}\subseteq \mathbb{D}_{2,\norm{\cdot}}$. If $\norm{\cdot}_@$ is a best ellipsoid norm for $\norm{\cdot}$ then $M_{\lVert \cdot \lVert}$ is the same as in (\ref{equation3}). If \,$\norm{(x,y)}=\sqrt{ax^2+cy^2}$ then $M_{\lVert \cdot \lVert}$ is best possible and we have $\mathbb{D}_{2,\norm{\cdot}}=[0,\gamma_{\norm{\cdot}}]$. In particular, the results in table (\ref{equation4}) hold. \qed

\end{Corollary}
\vspace*{0.01cm}
\begin{Corollary}
Let $\norm{\cdot}$ be an arbitrary norm on $\mathbb{R}^2$. Then the set of Dirichlet converging vectors (w.r.t $\norm{\cdot}$) is uncountable.\qed
\end{Corollary}

\ \\  We present another corollary of Theorem 1 which is formulated in the language of dynamics on the space of lattices.

\ \\ For $\norm{\cdot}$ an arbitrary norm on $\mathbb{R}^d$, we define a norm $\norm{\cdot}^*$ on $\mathbb{R}^{d+1}$ by\\ \,$\boldsymbol{\norm{(\vec{x},z)}^*}:=\max\{\,\norm{\vec{x}},|z|\,\}$.

\ \\ Let $\mathbb{X}_{d+1}$ be the space of lattices in $\mathbb{R}^{d+1}$ of covolume 1, and for a lattice $L\in\mathbb{X}_{d+1}$ we define \textbf{Minkowski's first successive minima} (w.r.t $\norm{\cdot}^*$) to be\\ $\boldsymbol{\lambda_{1}(L,\norm{\cdot}^*)}:=\min\{\,\norm{\vec{u}}^*\,|\,\vec{0}\neq\vec{u}\in L\}$, which is well defined by discreteness.

\ \\ For $\varepsilon>0$ define 
$\boldsymbol{\mathbb{K}_{\varepsilon}^{\norm{\cdot}^*}}:=\{L\in\mathbb{X}_{d+1}\,|\,\varepsilon\leq\lambda_{1}(L,\norm{\cdot}^*)\,\}$. \ \\ In particular, by Mahler's compactness criterion we have that $\{\mathbb{K}_{\varepsilon}^{\norm{\cdot}^*}\}_{\varepsilon>0}$ is an exhaustion by compact sets of $\mathbb{X}_{d+1}$. \ \\ For $t\geq0$, we denote by $g_t$ the diagonal $(d+1)\times(d+1)$ matrix\\\, $\boldsymbol{g_t}:=\,$diag$\,(e^t,\,...\,,e^t,e^{-dt})$.

\ \\ From a dynamical point of view, Theorem 1 allows us to construct 3-dimensional lattices such that the trajectory $g_tL$ has, in some sense, prescribed geometry. For example, it allows us to deduce the following:

\begin{Corollary}
  Let $\norm{\cdot}$ be an arbitrary norm on $\mathbb{R}^2$. \\Then for $M_{\norm{\cdot}}$ as in (\ref{equation3}), for all $r\in[0,M_{\norm{\cdot}}^{\frac{2}{3}}]$ there exist continuum many lattices $L\in\mathbb{X}_3$ such that
$$r=\min\,\{ \,\varepsilon>0\,\,|\,\,\,g_tL\notin\mathbb{K}_{\varepsilon}^{\norm{\cdot}^*}\,for\,\,all\,\, t\,\,large\,\,enough\}.
$$  
\end{Corollary} \textit{Proof.} \,For $\vec{v}=(v_1,v_2)\in\mathbb{R}^2$ define $l_{\vec{v}}:=\begin{pmatrix}
1 & 0 & -v_1\\
0 & 1 & -v_2\\
0 & 0 & 1
\end{pmatrix},\,L_{\vec{v}}:=l_{\vec{v}}\,\mathbb{Z}^3\in\mathbb{X}_3.$

\ \\Let $r\in[0,M_{\norm{\cdot}}^{\frac{2}{3}}]$, and by Theorem 1 let $\vec{v}\in\mathbb{R}^2$ be one of the uncountably many vectors such that for all $n\geq0$ we have $\sqrt{q_{n+1}}\,\lVert q_n\vec{v}-\vec{p}_n\lVert\in[r-\frac{1}{n+1},r-\frac{1}{n+2}]^{\frac{3}{2}}$.

\ \\ Since the map $t\longmapsto\lambda_{1}(g_tL_{\vec{v}}\,,\norm{\cdot}^*)$ is continuous and $\mathbb{Z}^3$ is discrete, there exist a sequence $(t_n,\vec{z}_n,Q_n)_{n=0}^{\infty}\subset\mathbb{R}_{\geq0}\times\mathbb{Z}^2\times\mathbb{N}$ with $0=t_0<t_1<...<t_n \rightarrow\infty$ such that if $t\in[t_{n},t_{n+1}]$  then $\lambda_{1}(g_tL_{\vec{v}}\,,\norm{\cdot}^*)=\lVert g_t\,l_{\vec{v}}\,(\vec{z}_n,Q_n) 
\lVert^*$, where the vectors $(\vec{z}_n)_{n=0}^{\infty}$ are not uniquely defined only up to (at most) a finite number of indices.

\ \\ One can then show that up to a finite number of indices, in fact we have that $(\vec{z}_n,Q_n)=(\vec{p}_n,q_n)$ -- i.e. the sequence of vectors which realises $\lambda_{1}(g_tL_{\vec{v}}\,,\norm{\cdot}^*)$ is exactly the sequence of best approximations of $\vec{v}$
(see \cite{Che} chapter 2). 

\ \\ It turns out that the sequence $(e^{t_{n+1}}\lVert q_{n}\vec{v}-\vec{p}_{n}\lVert)_{n=0}^{\infty}=(e^{-2t_{n+1}}q_{n+1})_{n=0}^{\infty}$ is the sequence of all local maxima of $\big(\lambda_{1}(g_tL_{\vec{v}}\,,\norm{\cdot}^*)\big)_{t>0}$, with the exception that $\lVert q_0\vec{v}-\vec{p}_0\lVert=\lVert \vec{v}-\vec{p}_0\lVert$ might also be a maximum at the beginning of the trajectory which correspond to $\lambda_{1}(g_0L_{\vec{v}}\,,\norm{\cdot}^*)$ in the case that $\lVert \vec{v}-\vec{p}_0\lVert>1$.

\ \\ In particular for all $n\geq0$ we have 
\begin{align*}
 \lambda_{1}(g_{t_{n+1}}L_{\vec{v}}\,,\norm{\cdot}^*)^3=q_{n+1}\,\lVert q_{n}\vec{v}-\vec{p}_{n}\lVert\,^2\in[r-\tfrac{1}{n+1},r-\tfrac{1}{n+2}]^3.   
\end{align*} So for all $t$ with $t\geq t_1$ we have $\lambda_{1}(g_{t}L_{\vec{v}}\,,\norm{\cdot}^*)<r$, and for all $n\geq0$ we have $\lambda_{1}(g_{t_{n+1}}L_{\vec{v}}\,,\norm{\cdot}^*)\in[r-\tfrac{1}{n+1},r-\tfrac{1}{n+2}]$. This finishes the proof. 
\qed 

$\vspace*{0.3cm}$ \ \\ 
\noindent \textbf{Remark 4.} One should notice that $\frac{2}{3}=\frac{d}{d+1}$ for $d=2$. 

\noindent The proof from above can be easily generalised to higher dimensions, and it shows that for $\norm{\cdot}$ an arbitrary norm on $\mathbb{R}^d$, controlling the lengths of $\vec{\beta}_n(\vec{v})$ is the same as controlling the upper bounds of the values of $\lambda_{1}(g_{t}L_{\vec{v}}\,,\norm{\cdot}^*)^{\frac{d}{d+1}}$.

\section{Reducing Theorem 1 to Theorem 2 }
\ \\ 
For the sake of convenience, from now on we denote the best approximations sequence of $\vec{v}$ by $((q_n),(\vec p_n)) $, omitting obvious indices. 
\vspace*{0.2cm}
\begin{Proposition}\label{Proposition1}
Let $\lVert \cdot \lVert$ be a norm on $\mathbb{R}^d$, and define a new norm $\boldsymbol{\lVert \cdot \lVert_\lambda}:=\lambda\,\lVert \cdot \lVert$ for some positive constant $\lambda>0$. Let $\vec{v}\in \mathbb{R}^d$. Then the following are equivalent:  
\begin{enumerate}
    \item $\big((q_n),(\vec p_n)\big)$ is the best approximations sequence of $\vec{v}$ w.r.t to $\lVert \cdot \lVert$.
    \item $\big((q_n),(\vec p_n)\big)$ is the best approximations sequence of $\vec{v}$ w.r.t to $\lVert \cdot \lVert_{\lambda}$.
\end{enumerate}
Furthermore, if $\norm{\cdot}_@$ is a best ellipsoid norm for $\norm{\cdot}$ with parameters $a,c,D$, then $\norm{(u_1,u_2)}_{\lambda,@}:=\sqrt{a\lambda^2u_1^2+c\lambda^2u_2^2}$ is a best ellipsoid norm for $\lVert \cdot \lVert_\lambda$ with parameters $a\lambda^2,c\lambda^2,D$.

\end{Proposition}

\vspace*{0.2cm}\begin{Proposition}
Let $\lVert \cdot \lVert'$ be a  norm on $\mathbb{R}^2$ such that there exists a best ellipsoid norm $\lVert \cdot \lVert_@'$ with parameters $a',c',D'$ (with respect to $\norm{\cdot}'$) with the additional property that $a'=1$. \\ If Theorem 1 holds for $\norm{\cdot}'$ then it holds for the general case.
\end{Proposition}

\ \\ We leave the proofs of Propositions 1 and 2 as exercises for the reader.

$\vspace*{0.3cm}$\ \\ 
\noindent Given $\vec{v} \in \mathbb{R}^{2} \setminus \mathbb{Q}^{2}$, a norm $\lVert \cdot \lVert$ and $\big((q_n),(\vec p_n)\big)_{n=0}^{\infty}$ the best approximations sequence of $\vec{v}$ (w.r.t to $\lVert \cdot \lVert$), define for all $1\leq n$
$$ \boldsymbol{\Pi_{n}(\vec{v})} :=  \{(\alpha,\vec{y}) \in\mathbb{R}^{3} \mid 0 \leq \alpha \leq q_{n} \,\,\, and \,\,\, \lVert \alpha\vec{v}-\vec{y}\, \lVert \, \leq \lVert q_{n-1}\vec{v}-\vec{p}_{n-1} \lVert \, \}.$$ 

\vspace*{0.1cm}
\begin{Proposition}\label{Proposition4}
    Let $\vec{v}\in\mathbb{R}^2$ and let $(q_n,\vec p_n)_{n=0}^{\infty}\subset \mathbb{N} \times \mathbb{Z}^{2}$ be a sequence such that the followings hold:
    
    \item \begin{itemize}
        \item $q_0=1$
        \item $q_n<q_{n+1}$ \,\,\,\,\,\,\,\,\,\,\,\,\,\,\,\,\,\,\,\,\,\,\,\,\,\,\,\,\,\,\,\,\,\,\,\,\,\,\,\,\,\,\,\,\,\,\,\,\,\,\,\,\,\,\,\,\,\,\, 
        \item $\lVert q_{n}\vec{v}-\vec{p}_{n} \lVert \, < \lVert q_{n-1}\vec{v}-\vec{p}_{n-1} \lVert $ 
        \item $int \, \Pi_{n}(\vec{v}) \cap \mathbb{Z}^{3}= \varnothing. $
    \end{itemize}

\ \\ Then $(q_n,\vec p_n)_{n=0}^{\infty}$ is the best approximations sequence of $\vec{v}$.
\end{Proposition}

\ \\ \textit{Proof.}\,\,
 Assume by contradiction there exist $Q\in\mathbb{N}$ with $q_{n-1}<Q<q_n$ for some $n$ and with $ \langle Q\vec{v}\rangle_{\lVert \cdot \lVert}<\langle q_{n-1}\vec{v}\rangle_{\lVert \cdot \lVert} $. So there exist $\vec{P}\in\mathbb{Z}^2$ such that \\$ \lVert Q\vec{v}-\vec{P} \lVert < \lVert q_{n-1}\vec{v}-\vec{p}_{n-1} \lVert$. So $(Q,\vec{P})\in int \, \Pi_{n}(\vec{v}) \cap \mathbb{Z}^{3}$, contradicting the assumptions. \\ Assume by contradiction $\vec{p}_n$ does not realise $\langle q_n\vec{v}\rangle_{\lVert \cdot \lVert}$. So there exist $\vec{s}\in\mathbb{Z}^2$ such that $\lVert q_{n}\vec{v}-\vec{p}_{n} \lVert \, > \lVert q_{n}\vec{v}-\vec{s}\, \lVert $. So $(q_n,\vec{s})\in int \, \Pi_{n+1}(\vec{v}) \cap \mathbb{Z}^{3} $ contradicting the assumptions.\ \\
\qed 

\ \\ We will deduce Theorem 1 from Theorem 2, to be stated below, which will be proved by an inductive procedure. In detail, we will construct by induction a sequence
$$\boldsymbol{\vec{w}_n} :=(q_n, \vec{p}_n)_{n=0}^{\infty} \in \mathbb{N} \times \mathbb{Z}^{2} $$
satisfying various properties which will imply that this sequence is the best approximation vectors sequence of the vector $\vec{v} \in \mathbb{R}^2$ defined by

\begin{equation}\label{equation5}
\boldsymbol{\vec{v}_n}:=  \,\frac{1}{q_n} \, \vec{p}_n  \,\,\hspace*{1cm}
\boldsymbol{\vec{v}} :=  \lim_{n\to\infty}
\vec{v}_n
\end{equation} and where our construction will imply that the limit exists. These properties will further imply that $\vec{v}$ and its best approximations sequence satisfy properties A,B,C from Theorem 1. Before doing so, we shall introduce some more definitions. 

\ \\Let $\norm{\cdot}_@$ 
 be a best ellipsoid norm for $\norm{\cdot}$ with parameters $a,c,D$. For $\vec{u} \in \mathbb{R}^2 $\, and positive real numbers $E,R$ denote
\begin{equation*}
\boldsymbol{ \Pi\,\biggl(\vec{u}\,,\,E\,,\,R\biggr)} :=  \left\{(\alpha,\vec{y}) \in\mathbb{R}^{3} \mid 0 \leq \alpha \leq E \,\,\, and \,\,\, \lVert \alpha\vec{u}-\vec{y}\, \lVert \, \leq R \, \right\}     
\end{equation*} \begin{equation*} 
\boldsymbol{ \Pi_{@}\,\biggl(\vec{u}\,,\,E\,,\,R\biggr)} :=  \left\{ (\alpha,\vec{y}) \in\mathbb{R}^{3} \mid 0 \leq \alpha \leq E \,\,\, and \,\,\, \lVert \alpha\vec{u}-\vec{y}\, \lVert_{@} \, \leq R \, \right\}.
\end{equation*}

\newpage \noindent Using this, and assuming \, $\vec{w}_{0},\vec{w}_{1},...,\vec{w}_{n-1}, \vec{w}_{n}$ \, have already been constructed, $\vec{v}_n$ is the same as in (\ref{equation5}), for $j$ $\in$ $\{1,...,n\}$ define further

$\hspace*{0.4cm}$  

$\boldsymbol{\vec{\psi}_{n,j}}:=q_{j-1}\vec{v}_{n}-\vec{p}_{j-1} \hspace*{1.5cm}
\boldsymbol{R_{n,j}}:= \lVert \vec{\psi}_{n,j}     \lVert$

$\hspace*{1.8cm}$

$\boldsymbol{R_{n,j,@}}:= \lVert \vec{\psi}_{n,j}     \lVert_@\hspace*{2.1cm}\boldsymbol{\Pi_{n,j} }:= \Pi\,(\vec{v}_{n},q_{j},R_{n,j})$.

\vspace*{0.3cm}
\begin{Theorem}
\begin{flushleft}
\ \\    
\ \\ 
Let $\lVert \cdot \lVert$ be a  norm on $\mathbb{R}^2$ such that there exists a best ellipsoid norm $\lVert \cdot \lVert_@$ with parameters $a,c,D$ with the additional property that $a=1$. Then there exists a positive constant $M_{\lVert \cdot \lVert}>0$ such that the following holds: 

\ \\ 
Let $\theta$,$\delta$ with $0\leq\theta<2\pi$ and $\,0<\delta$.

\ \\ Let $(I_{n})_{n=0}^{\infty}$ be a sequence of open intervals with $I_{n}\subset[0,M_{\lVert \cdot \lVert}]$ for all n.

\ \\ Let $m\in\mathbb{N}$, and let $(z_{n})_{n=1}^{\infty}\in\{0,1,...,m-1\}$ be a sequence of residues mod $m$ such that $z_n\in\{sz_{n-1}+iz_{n-2}\pm z_{n-3}\mod{m}\,|\,\,i,s\in\mathbb{Z}\}$ for all $n$, where $z_0=1$, $z_{-1}=z_{-2}=0$.

\vspace*{0.4cm} 
\end{flushleft}

Then there exist continuum many sequences $\vec{w}_n =(q_n, \vec{p}_n)_{n=0}^{\infty} \in \mathbb{N} \times \mathbb{Z}^{2}$ satisfying the following properties:

\begin{enumerate}[(2.1)]
    \item 
    \begin{itemize}
        \item $\Pi_{n,n} \cap \mathbb{Z}^3 = 
    \{\vec{w}_n, \vec{w}_{n-1}, \vec{w}_n-\vec{w}_{n-1},\vec{0} \}$
    $\hspace*{0.57cm}$ $1\leq n$ 
    \item $\Pi_{n,j} \cap \mathbb{Z}^3 \subseteq 
    \{\vec{w}_j, \vec{w}_{j-1},\vec{w}_j-\vec{w}_{j-1},\vec{0}\}$
\,\,\,\,\,\,\,\,\,\,\,\,\, $1\leq j \leq n-1$
    \item   $int\,\Pi_{n,j} \cap \mathbb{Z}^3 =\varnothing$ \hspace*{3.74cm} $1\leq n$ and\,\,$1\leq j \leq n$
    \end{itemize}
    
    \label{(2.1)}
  
    \item $n q_{n-1}<q_n$ \hspace*{4.75cm} for all $1\leq n$
    $\vspace*{0.15cm}$ \\$q_0=1$
    
    \label{(2.2)}
   
    \item  $2R_{n,j}<R_{n,j-1}$ \hspace*{4.1cm} for all $1\leq n$ and\, $1\leq j \leq n$
    \label{(2.3)}
    
    \item $ \sqrt{q_{j}}\,\,   \lVert q_{j-1}\vec{v}_n-\vec{p}_{j-1} \lVert\,\, \in \,int\,I_{j-1}$ \hspace*{1.5cm} for all $1\leq n$ and\, $1\leq j \leq n$
    \label{(2.4)}
    
    \item $\lVert \vec{v}_n-\vec{v}_{n-1}\lVert < 2^{-n}$ \hspace*{3.55cm} for all $1\leq n$
    \label{(2.5)}
    
    \item $\vec{\psi}_{n,j}\in \,int\, \Omega(\theta,\delta)$ \hspace*{3.68cm} for all $1\leq n$ and\, $1\leq j \leq n$
    $\vspace*{0.15cm}$ \\ $\vec{\psi}_{n,j}$ and $\vec{\psi}_{n,j+1}$ always lie in opposite quadrants
    \label{(2.6)}
    
    \item $q_n\equiv z_{n} \, mod\,\, m$  \hspace*{3.95cm} for all $1\leq n$.  
    \label{(2.7)}

\end{enumerate} 
Furthermore, $\{\vec{w}_{n-2},\vec{w}_{n-1},\vec{w}_n\}$ is a basis of the lattice $\mathbb{Z}^3$ for all $n\geq2$, and the above holds for 
$$ M_{\lVert \cdot \lVert}= \big(\frac{4c}{4D^2-1}\big)^{\frac{1}{4}}.$$
\end{Theorem}  

\ \\ \textit{Proof of Theorem 1 assuming Theorem 2}. \,\,By Proposition 3 it's enough to prove Theorem 1 for the special case that $\norm{\cdot}$ has a best ellipsoid norm $\lVert \cdot \lVert_@$ with parameters $a,c,D$ with the additional property that $a=1$. \\ By Theorem 2, let $\vec{w}_n =(q_n, \vec{p}_n)$ be a sequence with properties (2.1)-(2.7). 

\ \\ Define ${\vec{v}} := \lim_{n\to\infty}\vec{v}_n$ which is well defined by property (2.5). \ \\ First we show: 
\begin{enumerate}[(a)]
    \item $(q_n, \vec{p}_n)$ is a sequence of best approximations of $\vec{v}$.
    \item properties (A),(B),(C) hold.
\end{enumerate}
To prove (a) we use Proposition \ref{Proposition4}:
\begin{itemize}
    \item $q_n<q_{n+1}$ by property (2.2).
    \item For $1\leq j\in \mathbb{N}$ we have
    
\begin{dontbotheriftheequationisoverlong}
 \lVert q_{j}\vec{v}-\vec{p}_{j} \lVert = 
\lVert q_{j}\lim_{n\to\infty}\vec{v}_n-\vec{p}_{j} \lVert = 
\lim_{n\to\infty}\lVert q_{j}\vec{v}_n-\vec{p}_{j} \lVert \,<
\lim_{n\to\infty}\lVert q_{j-1}\vec{v}_n-\vec{p}_{j-1} \lVert = 
\lVert q_{j-1}\vec{v}-\vec{p}_{j-1} \lVert   
\end{dontbotheriftheequationisoverlong}

where the inequality follows from property (2.3) and the fact that both sequences converge to a positive number.

     \item Assume by contradiction that there exist $1\leq j\in \mathbb{N}$ and $(Q,\vec{P})\in \mathbb{N} \times \mathbb{Z}^{2}$ such that $(Q,\vec{P})\in int \, \Pi_{j}\cap \mathbb{Z}^{3}$. I.e. we have
     
     $$1\leq Q < q_{j}  \,\, and \,\,\,
     \lVert Q\vec{v}-\vec{P} \lVert <
     \lVert q_{j-1}\vec{v}-\vec{p}_{j-1} \lVert.$$
     
     In particular, for $n$ which is large enough,
     $$ 1\leq Q < q_{j}  \,\, and \,\,\,
     \lVert Q\vec{v}_n-\vec{P} \lVert <
     \lVert q_{j-1}\vec{v}_n-\vec{p}_{j-1} \lVert$$
     meaning that $(Q,\vec{P}) \in int \, \Pi_{n,j}\cap \mathbb{Z}^{3}$, contradicting property (2.1). 
\end{itemize} We deduce by Proposition \ref{Proposition4} that $(q_n, \vec{p}_n)$ is the sequence of best approximations of $\vec{v}$.

\ \\ For the proof of properties (A), (B), (C) we have
\begin{itemize}
    \item \,\,\,\, 
    $ \sqrt{q_{j}}\,\,   \lVert q_{j-1}\vec{v}-\vec{p}_{j-1} \lVert\,\, = \lim_{n\to\infty}\sqrt{q_{j}}\,\,   \lVert q_{j-1}\vec{v}_n-\vec{p}_{j-1} \lVert.\,\,
    $

    According to (2.4) the RHS is inside $int\, I_{j-1}$ for all $n>j$. As $I_{j-1}$ is a closed interval, we get (A).
    \item\,\,\,\,
    $\vec{\beta}_j(\vec{v})=\lim_{n\to\infty}\sqrt{q_{j+1}}\,\vec{\psi}_{n,j+1}$. 
    
    \ \\ According to (2.6) $\vec{\psi}_{n,j+1}\in int\, \Omega(\theta,\delta)$ for all $j\leq n-1$. $\Omega(\theta,\delta)$ is close, so we have $\vec{\beta}_j(\vec{v})\in \Omega(\theta,\delta)$. 
    Furthermore, as $\vec{\psi}_{n,j}$ and $\vec{\psi}_{n,j+1}$ always lie in opposite quadrants, we also have that $\vec{\beta}_{j-1}(\vec{v})$ and $\vec{\beta}_j(\vec{v})$ lie in opposite quadrants, so we get (B). \item \,\,\,\, (C) follows immediately from (2.7).
    
\end{itemize}
\ \\ As\, $\sup_{n\in\mathbb{N}
}\dfrac{q_{n+1}}{q_n}=\infty$ we have that $\vec{v}$ is not badly approximable (see \cite{AM2}).

\ \\ If two vectors have different sequences of approximation vectors, even after possibly altering finitely many elements in the sequences, then they are different. So from the continuum many sequences we get in Theorem 2 we get continuum many vectors which satisfy the conclusions of Theorem 1. This finishes the proof of Theorem 1. \,\,\qed

\section{Proving Theorem 2 -- simplifying the settings}
For the base of the induction we set
$$ \vec{w}_{0}:=(1,0,0)\,\,\,\,\,\vec{w}_{-1}:=(0,p_{-1,1},p_{-1,2}) \,\,\,\,\,\vec{w}_{-2}:=(0,p_{-2,1},p_{-2,2}) $$
where we choose $(p_{-1,1},p_{-1,2})\in\mathbb{Z}^2$ to be two co-prime integers such that $(p_{-1,1},p_{-1,2})\in int\, \Omega(\theta,\delta)$, and we choose $(p_{-2,1},p_{-2,2})$ to be two integers such that $\vec{w}_{-2}$ completes $\{\vec{w}_0,\vec{w}_{-1}\}$ to a basis of $\mathbb{Z}^3$. This choice is needed in order to construct $\vec{w}_1,\vec{w}_2$ and $\vec{w}_3$, and will become clear throughout the paper.

\ \\ Assume by induction that we have $\vec{w}_0,\vec{w}_1,...\,,\vec{w}_{n-2},\vec{w}_{n-1}$ which satisfy properties (2.1) to (2.7). 

\ \\ Given two linearly independent integer vectors $\vec{u}_1,\vec{u}_2 \in\mathbb{Z}^3$, we say that $\{\vec{u}_1,\vec{u}_2\}$ \textbf{is a primitive set} if there exist $\vec{u}_3$ such that $\{\vec{u}_1,\vec{u}_2,\vec{u}_3\}$ is a basis of the lattice $\mathbb{Z}^3$. Equivalently, if $span_{\mathbb{R}}\{\vec{u}_1,\vec{u}_2\}\cap\mathbb{Z}^3=span_{\mathbb{Z}}\{\vec{u}_1,\vec{u}_2\}$.

\ \\ Now define
$$ \boldsymbol{V}:=span_{\mathbb{R}}\,\{\vec{w}_{n-1},\vec{w}_{n-2}\} $$
$$
\boldsymbol{\Gamma}:=V\cap\mathbb{Z}^3=span_{\mathbb{Z}}\,\{\vec{w}_{n-1},\vec{w}_{n-2}\}
$$ where the last equality follows from the inductive assumption that \ \\ $\{\vec{w}_{n-3},\vec{w}_{n-2},\vec{w}_{n-1}\}$ is a basis of the lattice $\mathbb{Z}^3$, hence $\{\vec{w}_{n-2},\vec{w}_{n-1}\}$ is a primitive set. 

\ \\Without loss of generality, let $i_0,s_0\in\{0,1,...,m-1\}$ such that 

\begin{equation}\label{equation666}
z_n\equiv s_0z_{n-1}+i_0z_{n-2}+z_{n-3}\mod{m}  
\end{equation} \ \\ and define further
$$ \boldsymbol{V'}:=V+\vec{w}_{n-3} $$
$$ \boldsymbol{\Gamma'}:=V'\cap\mathbb{Z}^3
=\Gamma+\vec{w}_{n-3}.$$ In case that $z_n\equiv s_0z_{n-1}+i_0z_{n-2}\boldsymbol{-}z_{n-3}\mod{m}$, define $V':=V-\vec{w}_{n-3}$ and $\Gamma':=\Gamma-\vec{w}_{n-3}$. This choice of sign is only needed in order to prove (2.\ref{(2.7)}), the rest of the proof remains identical in both cases.

\newpage \ \\ Let $\norm{\cdot}_@$
be a best ellipsoid norm (w.r.t $\norm{\cdot}$)
with parameters $a,c,D$ with the additional property that $a=1$. Let $\varphi$ be the projection from $\mathbb{R}^3$ on the second two coordinates (i.e. $\varphi$ is represented by $\begin{pmatrix}
0 & 1 & 0\\
0 & 0 & 1
\end{pmatrix}$\,). 

\ \\ Let $\vec{w}_n=(q_n,\vec{p}_n)\in\mathbb{N}\times\mathbb{Z}^2$ be any future choice of $\vec{w}_n$. For $G\in GL_3(\mathbb{R})$ and $j\in\{0,1...,n\}$ define

$$\boldsymbol{\widetilde{p_j}} :=\varphi\,\bigl( G\bigl( \vec{w}_j \bigr) \bigr)\,\,\,\,\,\,\,\,\,
\boldsymbol{\widetilde{v_j}} :=\varphi\,\bigl( G\bigl( 1,\vec{v}_j \bigr) \bigr)=\frac{\widetilde{p_j}}{q_j}\,\,\,\,\,\,\,\,\,\boldsymbol{\widetilde{w}_j}:=(q_j,\widetilde{p}_j).$$

\begin{Lemma}
\label{Proposition5}
Let $\vec{w}_0,\vec{w}_1,...\,,\vec{w}_{n-2},\vec{w}_{n-1}$ be the integer vectors which satisfy properties (2.1) to (2.7). Then there exists $G\in GL_3(\mathbb{R})$ and two positive numbers $H,L$ such that for any future choice of $\vec{w}_n=(q_n,\vec{p}_n)\in\mathbb{N}\times\mathbb{Z}^2$ the properties below hold:

\begin{enumerate}[(G.1)]
\item
$G=\begin{pmatrix}
1 & 0 & 0\\
* & * & *\\
* & * & *
\end{pmatrix}$.
\item
The lower-right minor of G is orthogonal with respect to $\norm{\cdot}_{@}$.\\ 
I.e. if $G=\begin{pmatrix}
1 & 0 & 0\\
* & g_1 & g_2\\
* & g_3 & g_4
\end{pmatrix}$ then $\norm{g(\vec{u})}_@=\norm{\vec{u}}_@$ for   
$g=\begin{pmatrix}
g_1 & g_2\\
g_3 & g_4 
\end{pmatrix}$ and for all $\vec{u}\in\mathbb{R}^2$.

\item
$G(\vec{w}_{n-1})=(q_{n-1},0,0)$.
\item
$G(V)\,=\,\,\{\,(u_1,u_2,0)\,\mid\,(u_1,u_2)\in\mathbb{R}^2
\}.
$
\item
$G(V')=\,\{\,(u_1,u_2,H)\mid\,(u_1,u_2)\in\mathbb{R}^2
\}.$ 

\item
$G(\vec{w}_{n-2})=(q_{n-2},L,0)$.

\item $G(V)$ and $G(V')$ are neighbouring planes with respect to the lattice $G(\mathbb{Z}^3)$. I.e. $G(V')=G(V)+G(\vec{w}_{n-3})$ and for all $0<t<1$ we have that \\ $\bigl(G(V)+tG(\vec{w}_{n-3})\bigr)\cap G(\mathbb{Z}^3)=\varnothing $.

\item
    $G(\Gamma)=span_{\mathbb{Z}}\{G(\vec{w}_{n-1}),G(\vec{w}_{n-2})\}=span_{\mathbb{Z}}\{(q_{n-1},0,0),(q_{n-2},L,0)\}.$
    $\vspace*{0.15cm}$ \\ Additionally, define $Covol(G(\Gamma))$ to be the area of the parallelogram \\ $\{\alpha G(\vec{w}_{n-1})+\beta G(\vec{w}_{n-2})\mid (\alpha,\beta)\in[0,1)^2\}$. $\vspace*{0.15cm}$ \\ Then
    $Covol(G(\Gamma))=q_{n-1}L=H^{-1}$.
    
\item $G(\Gamma ')=G(\Gamma + \vec{w}_{n-3}) = G(\Gamma)+\widetilde{w}_{n-3}$.
    
\item 
    If $\vec{U}=(\alpha,\vec{u})\,,\vec{T}=(\alpha,\vec{t})$ then we have that\\ 
    $\norm{\vec{u}-\vec{t}\,}_@=\lVert\varphi(G(\vec{U}))-\varphi(G(\Vec{T}))\lVert_@=\lVert\varphi(G^{-1}(\vec{U}))-\varphi(G^{-1}(\Vec{T}))\lVert_@$. 
    
    In particular, 
    
    for $j\in\{0,1,...,n\}$ we have $\lVert \, \vec{v}_{n}-\vec{v}_{j-1}   \lVert_{@} \,\, = \, \lVert \, \widetilde{v}_{n}-\widetilde{v}_{j-1}\,  \lVert_{@}$, 
    
    \hspace*{0.3cm} \, hence  $R_{n,j,@}=\lVert \, q_{j-1}\vec{v}_{n}-\vec{p}_{j-1}   \lVert_{@} \,\,= \lVert \, q_{j-1}\widetilde{v}_{n}-\widetilde{p}_{j-1} \,   \lVert_{@} $.

\item \label{G.11}$ 
\Pi_{n,j}\subseteq\Pi_{@}\,(\vec{v}_n,q_j,DR_{n,j,@})=G^{-1}\bigl(\Pi_@(\widetilde{v}_n,q_j,DR_{n,j,@})\bigr).$ 

In particular,\,\,$ \vec{u}\notin int\, \Pi_@(\,\widetilde{v}_n,q_j,DR_{n,j,@}) \Longrightarrow  G^{-1}(\vec{u})\notin int\,\Pi_{n,j}$. 
    
\end{enumerate}
 \end{Lemma}\ \\ \textit{Proof.}\,\,  First we show there exist $G\in GL_3(\mathbb{R})$ satisfying (G.1)-(G.6) and then we show (G.7)-(G.11) follow. Assume $G$ is of the form:
 $$G=\begin{pmatrix}
1 & 0 & 0\\
* & g_1 & g_2\\
* & g_3 & g_4
\end{pmatrix}=\begin{pmatrix}
1 & 0 & 0\\
0 & k_1 & k_2\\
0 & k_3 & k_4\\ 
\end{pmatrix}\begin{pmatrix}
1 & 0 & 0\\
0 & z_1 & z_2\\
0 & z_3 & z_4\\ 
\end{pmatrix} \begin{pmatrix}
1 & 0 & 0\\
x & 1 & 0\\
y & 0 & 1\\ 
\end{pmatrix}=KZA,$$ 
denote $
k=\begin{pmatrix}
k_1 & k_2\\
k_3 & k_4 
\end{pmatrix}$, $z=\begin{pmatrix}
z_1 & z_2\\
z_3 & z_4 
\end{pmatrix}$, and assume that $k,z\in GL_2(\mathbb{R})$.

\ \\ First choose the unique values of $x$ and $y$ such that $A(\Vec{w}_{n-1})=(q_{n-1},0,0)$. Then in particular, as $K,Z$ are of the form above, we now have (G.3) for any choice of $k,z\in GL_2(\mathbb{R})$.

\ \\ For (G.6), first notice that since $\vec{w}_{n-1}$ and $\vec{w}_{n-2}$ form a primitive pair, they are also linearly independent. So as $\det(A)=1$, we get that $A\vec{w}_{n-2}=(q_{n-2},*,*)\neq(q_{n-2},0,0)$. Now choose $z$ to be a rotation matrix that preserves $\norm{\cdot}_@$ (i.e. maps the ellipse which is the unit ball of 
$\norm{\cdot}_@$ back to itself) that satisfies $z\,\varphi \,A\,\vec{p}_{n-2}=(L,0)$ for some $0<L$. So we now have $ZA\vec{w}_{n-2}=(q_{n-2},L,0)$ and $ZA\vec{w}_{n-1}=(q_{n-1},0,0)$.

\ \\ As $G(V)=span_{\mathbb{R}}\{G(\vec{w}_{n-1}),G(\vec{w}_{n-2})\}=span_{\mathbb{R}}\{(q_{n-1},0,0),(q_{n-2},L,0)\}$ we immediately get (G.4). 

\ \\ As $\{\vec{w}_{n-1},\vec{w}_{n-2},\vec{w}_{n-3}\}$ is a basis of $\mathbb{R}^3$ and $G$ is invertible, $G(V')=\{\,(u_1,u_2,H)\\ \mid(u_1,u_2)\in\mathbb{R}^2
\}$ with $H$ a constant such that $0\neq H$. If $0<H$ we are done. if not, we choose the matrix $k$ to be a reflection matrix that preserves $\norm{\cdot}_@$ such that $k\,(0,L)=(0,L)$, and get that $G=KZA$ satisfies (G.1) to (G.6). 

\ \\We show properties (G.7)-(G.11) follow from (G.1)-(G.6) as follows: 
 \begin{itemize}
    \item (G.7) follows immediately from the definition of $V$, the fact that $\vec{w}_{n-3}$ completes the primitive set $\{\vec{w}_{n-1},\vec{w}_{n-2}\}$ into a basis of the lattice $\mathbb{Z}^3$, and the fact that G is a linear transformation.

    \item
     For (G.8), by definition of $Covol(G(\Gamma))$ we get that $Covol(G(\Gamma))=q_{n-1}L$. \\ To see why $q_{n-1}L=H^{-1}$, first notice that as $g$ is orthogonal with respect to $\norm{\cdot}_@$ we have that $$\begin{pmatrix}
a & 0\\
0 & c 
\end{pmatrix}=\begin{pmatrix}
g_1 & g_2\\
g_3 & g_4 
\end{pmatrix}^t\begin{pmatrix}
a & 0\\
0 & c 
\end{pmatrix}\begin{pmatrix}
g_1 & g_2\\
g_3 & g_4 
\end{pmatrix},$$
and so $\det(g)=\pm 1$. So we have \\   $Covol(G(\mathbb{Z}^3)=\bigl\lvert \det(G)\bigl\lvert\,=
     \biggl\lvert\,\det \begin{pmatrix}
1 & 0 & 0\\
* & g_1 & g_2\\
* & g_3 & g_4
\end{pmatrix}\biggl\lvert\,=1$. 
 
On the other side, $Covol(G(\mathbb{Z}^3))$ is also equal to the volume of the fundamental parallelpiped associated to a basis of $G(\mathbb{Z}^3)$. 

As 
$\{G(\vec{w}_{n-1}),G(\vec{w}_{n-2}),G(\vec{w}_{n-3})\}$ is a basis of $G(\mathbb{Z}^3)$, we get from (G.3),\\(G.5),(G.6) that
$\{(q_{n-1},0,0),(q_{n-2},L,0),(*,*,H)\}$ is a basis of $G(\mathbb{Z}^3)$ with corresponding fundamental parallelpiped having volume $q_{n-1}LH$.\\ So $q_{n-1}LH=1$.

 \item
    \vspace*{0.2cm} For (G.10),$$\lVert\varphi(G(\vec{U}))-\varphi(G(\Vec{T}))\lVert_@
\,=\,
\lVert\varphi(G(0,\vec{u}-\vec{t}))\lVert_@
\,=\,
\lVert g\,(\vec{u}-\vec{t})
     \lVert_{@}
\,\stackrel{(G.2)}{=}
\lVert\vec{u}-\vec{t}\lVert_@.
$$
\begin{dontbotheriftheequationisoverlong}
\lVert\varphi(G^{-1}(\vec{U}))-\varphi(G^{-1}(\Vec{T}))\lVert_@
\,=\,
\lVert\varphi(G^{-1}(0,\vec{u}-\vec{t}))\lVert_@
\,=\,
\lVert g^{-1}\,(\vec{u}-\vec{t})
     \lVert_{@}
\,\stackrel{(*)}{=}
\lVert\vec{u}-\vec{t}\lVert_@
\end{dontbotheriftheequationisoverlong} 

where for $(*)$ we use the fact that if $g$ if orthogonal with respect to $\norm{\cdot}_@$ then so is $g^{-1}$.
Applying the above to the vectors $(1,\vec{v}_n),(1,\vec{v}_{j-1})$ and multiplying the equation by $q_{j-1}$, we get (G.10).    
   
   \item
   \vspace*{0.2cm} For (G.11), $\Pi_{n,j}\subseteq\Pi_{@}\,(\vec{v}_n,q_j,DR_{n,j,@})$ is clear from the definitions and the fact that $\norm{\cdot}_@$ is a best ellipsoid norm for $\norm{\cdot}$. Additionally,
   \begin{flushleft}
    $\Pi_{@}\,(\vec{v}_n,q_j,DR_{n,j,@})$   
   \end{flushleft}
    \begin{flalign*}
    =\left\{ (\alpha,\vec{y}) \in\mathbb{R}^{3} \mid 0 \leq \alpha \leq q_j \,\,\, and \,\,\, \lVert \alpha\vec{v}_n-\vec{y}\, \lVert_@ \, \leq D\lVert \, q_{j-1}\vec{v}_{n}-\vec{p}_{j-1}   \lVert_@ \, \right\}
   \end{flalign*}
   \begin{equation}\label{equation6}
   \vspace*{1cm}
   \stackrel{(G.10)}{=}
   \left\{ (\alpha,\vec{y}) \in\mathbb{R}^{3} \mid 0 \leq \alpha \leq q_j \, and \, \lVert \alpha\vec{v}_n-\vec{y}\, \lVert_@ \, \leq D\lVert \, q_{j-1}\widetilde{v}_{n}-\widetilde{p}_{j-1}   \lVert_@ \, \right\}.    
   \end{equation} 
   $G$ is invertible, so there exist a unique $\vec{t}_y$ such that $G^{-1}(\alpha,\vec{t_y})=(\alpha,\vec{y})$. Hence the right hand side of (\ref{equation6}) is equal to
\begin{equation*}  
   \left\{   
G^{-1}(\alpha,\vec{t}_y)\in\mathbb{R}^{3}
 \,\,\middle\vert 
\begin{array}{l}
 \hspace*{3.2cm}0 \leq \alpha \leq q_j
 \vspace*{0.1cm} \\ 
\lVert \alpha\vec{v}_n-\varphi(G^{-1}(\alpha,\vec{t}_y))\, \lVert_@ \, \leq D\lVert \, q_{j-1}\widetilde{v}_{n}-\widetilde{p}_{j-1}   \lVert_@
\end{array}\right\}
\end{equation*}
\begin{dontbotheriftheequationisoverlong}
 =
G^{-1}\left\{   
(\alpha,\vec{t}_y) \in\mathbb{R}^{3}
 \middle\vert 
\begin{array}{l}
 \hspace*{5.1cm}0 \leq \alpha \leq q_j
 \vspace*{0.1cm} \\ 
\lVert \,\varphi(G^{-1}(\alpha,\alpha\widetilde{v}_n))-\varphi(G^{-1}(\alpha,\vec{t}_y))\, \lVert_@ \, \leq D\lVert \, q_{j-1}\widetilde{v}_{n}-\widetilde{p}_{j-1}   \lVert_@
\end{array}\right\}   
\end{dontbotheriftheequationisoverlong}
\begin{flalign*}
&=
G^{-1}\left\{\bigl   
(\alpha,\vec{t}_y) \in\mathbb{R}^{3}
 \middle\vert 
\begin{array}{l}
 \hspace*{1.4cm}0 \leq \alpha \leq q_j
 \vspace*{0.1cm} \\ 
\lVert \alpha\widetilde{v}_n-\vec{t}_y\, \lVert_@ \, \leq D\lVert \, q_{j-1}\widetilde{v}_{n}-\widetilde{p}_{j-1}   \lVert_@
\end{array}\right\}
\\ 
&=G^{-1}\bigl(\Pi_@(\widetilde{v}_n,q_j,DR_{n,j,@})\bigr).
  \end{flalign*}
\qed 
\end{itemize}

We now aim to use Lemma \ref{Proposition5} in order to proceed with the inductive construction – choosing the vector $\vec{w}_n$. Specifically, we will look for  $\widetilde{w}_n=(q_n,\widetilde{p}_n)\in G(\Gamma^{'})$ satisfying certain properties in order to eventually define\\ $\vec{w}_n := G^{-1}(\widetilde{w}_n)$.

 \section[chapter 4]{Proving Theorem 2 -- choosing the vector $\vec{w}_n$ by choosing the vector  $\widetilde{w}_n$}

Let $\vec{w}_0,\vec{w}_1,...\,,\vec{w}_{n-2},\vec{w}_{n-1}$ be the integer vectors which satisfy properties (2.1) to (2.7), and let $G\in GL_3(\mathbb{R})$ be a map satisfying Lemma \ref{Proposition5}.
 
\ \\ For $(x,y,H)\in G(V')$ define the cylinder
\begin{align*} 
&\boldsymbol{\Pi_{1,@}\bigl[x,y,H\big]}:=
\Pi_@\biggl((\frac{y}{x},\frac{H}{x})\,\,,\,\, |x|\,\,,\,\,\frac{Dq_{n-1}\,\lVert \, (y,H) \,   \lVert_{@}}{|x|}\, \biggr).
\end{align*} Notice that if we set $(q_n,\widetilde{p}_n):=(x,y,H)$ for some $(x,y)$ with $0<x$ to be chosen later, then  $\Pi_{1,@}\bigl[q_n,\widetilde{p}_n\big]=\Pi_@(\,\widetilde{v}_n,q_n,DR_{n,n,@})$. Indeed, 
\begin{align*}
&\Pi_{1,@}\big[q_n,\widetilde{p}_n\big]=\Pi_@\bigl(\widetilde{v}_n\,,q_n\,,\,\dfrac{Dq_{n-1}\,\lVert \, \widetilde{p}_n \,   \lVert_{@}}{q_n}\, \bigr)
\\ &=\,\Pi_@\bigl(\widetilde{v}_n\,,q_n\,,\,D\,\lVert \, q_{n-1}\widetilde{v}_n \,   \lVert_{@}\, \bigr)   \stackrel{*}{=}\Pi_@\bigl(\widetilde{v}_n\,,q_n\,,\,D\,\lVert \, q_{n-1}\widetilde{v}_n-\widetilde{p}_{n-1} \,   \lVert_{@}\, \bigr)
\end{align*}
\begin{equation}\label{equation7} 
    \vspace*{0.4cm} \\ \hspace*{1.2cm}\stackrel{**}{=}\Pi_@\bigl(\widetilde{v}_n\,,q_n\,,\,D\,\lVert \, q_{n-1}\vec{v}_n-\vec{p}_{n-1} \,   \lVert_{@}\, \bigr)=\Pi_@(\,\widetilde{v_n},q_n,DR_{n,n,@})
\end{equation} where $\stackrel{*}{=}$ holds since \,$\widetilde{p}_{n-1}=\vec{0}$, and $\stackrel{**}{=}$ holds by Lemma \ref{Proposition5} property (G.10).
 
\ \\ So if we know that $G(\vec{u})\notin int\, \Pi_@(\,\widetilde{v_n},q_n,DR_{n,n,@})$
for all $\vec{u}\in\mathbb{Z}^3$, then by property (G.11) we get that $\vec{u}\notin int\,\Pi_{n,n}$ -- i.e. we get useful information in order to prove property (2.1) of Theorem 2. Furthermore, by calculating the volume of $\Pi_{1,@}\bigl[q_n,\widetilde{p}_n\big]$ we get useful information regarding the value of \\$q_n\,\lVert q_{n-1}\vec{v}_{n}-\vec{p}_{n-1}\lVert^2$ -- i.e. we get useful information regarding property (2.4) of Theorem 2. Indeed, if we denote by $\mathbb{B}$ the area of the unit ball corresponding to the norm $\norm{\cdot}_@$, then we have \begin{align*}
\mathrm{Vol}\biggl(\Pi_{1,@}\big[q_n,\widetilde{p}_n\big]\biggr) =\mathrm{Vol}\biggl(\Pi_@\biggl(\widetilde{v}_n\,,q_n\,,\,D\,\lVert \, q_{n-1}\vec{v}_n-\vec{p}_{n-1} \,   \lVert_{@}\, \biggr)\biggr)
\end{align*}
\begin{equation}\label{equation8}
=\,q_n\,
\mathbb{B}\,D^2\norm{q_{n-1}\vec{v}_n-\vec{p}_{n-1}}^2_@ =\,\mathbb{B}\,(D\Lambda)^2\,q_n\norm{q_{n-1}\vec{v}_n-\vec{p}_{n-1}}^2    
\end{equation}

\ \\ where 
$\Lambda=\frac{\,\,\,\,\norm{q_{n-1}\vec{v}_n-\vec{p}_{n-1}}_@}{\norm{q_{n-1}\vec{v}_n-\vec{p}_{n-1}}}$ with $\frac{1}{D}\leq \Lambda\leq1$. 

\ \\ In order to study cylinders of the form $\Pi_{1,@}\bigl[x,y,H\big]$ we introduce another parameterization via points from $G(V)$, and we introduce functions which allow us to go from one parameterization to another.

\newpage \ \\ For the sake of convenience denote from now on $\boldsymbol{q}:=q_{n-1}$,\\ $\boldsymbol{\Theta_{x}}:=\lVert(x,q)\lVert_@$, and define: 

\begin{align*} 
&\boldsymbol{\Pi_{2,@}\bigl[x,y,0\,\bigr]}\,\,:=
\Pi_@\biggl(\,(\frac{x\,y}{\Theta_{x}^2},\frac{q\,y}{\Theta_{x}^{2}})\,\,,\,\, \frac{H\,\Theta_{x}^{2}}{q\,|y|}\,\,,\,\,\frac{Dq\,\,|y|}{\Theta_{x}} \biggr)
\\ \\
&\boldsymbol{f_2\bigl(x,y,0\bigr)}:=
\bigl(\,\dfrac{H\,\Theta_{x}^2}{q\,y}\,,\,\frac{H\,x}{q}\,,\,H\,\bigr)
\,\,\,\,\,\,\,\,\,\,\,\,\,\,\,\,\,\,\,\,\,\,\, f_2:G(V)\longrightarrow G(V')
\end{align*} 

\ \\ The next lemma shows the connection between the cylinders $\Pi_{1,@}\bigl[x,y,H\big],$\\$\Pi_{2,@}\bigl[x,y,0\bigr]$ and the function $f_2$.
\begin{Lemma}
\label{Lemma2}
For every $(x,y,0)\in G(V)$ we have
$$\Pi_{1,@}\bigl[f_2(x,y,0)\bigr]= \Pi_{2,@}\big[x,y,0\big].$$
\end{Lemma}
 
\ \\ \textit{Proof.} \,\,The proof is done by direct calculation and is left as an exercise for the reader. \qed

\vspace*{0.3cm} \begin{Lemma}\label{Lemma3}
Let $(x,y,0)\in G(V)$ and let $(x',y',H)=f_2\bigl(x,y,0\bigr)\in G(V')$. 

\noindent Then $x'=\dfrac{q\lVert(y',H)\lVert_@^2}{yH}$.\, Furthermore,

\begin{align*}
Vol\biggl(\Pi_{1,@}\big[x',y',H\big]\biggr)=Vol\biggl(\Pi_{2,@}\big[x,y,0\big]\biggr)=D^2Hq\,\mathbb{B}\,|y|=\frac{D^2\,\mathbb{B}\,|y|}{L}.
\end{align*} \end{Lemma}

\vspace*{0.3cm}
\noindent \textit{Proof.} \,\,Proof is done again by simple calculations. 

\ \\ From the definition of $f_2(x,y,0)$ we have that $y'=\frac{Hx}{q}$ and $x'=\frac{H\,\Theta_{x}^2}{qy}$,\\ so we get that \begin{align*}
\frac{q\lVert(y',H)\lVert_@^2}{yH}=\frac{q\lVert(\,\dfrac{H\,x}{q},H)\lVert_@^2}{yH}=
\frac{\,\lVert(\,Hx,Hq)\lVert_@^2}{qyH}=
\frac{H\,\Theta_{x}^2}{qy}=x'.
\end{align*}

\ \\  For the volume calculation, by Lemma \ref{Lemma2} we can calculate the volume of the cylinders via each one of the parameterizations. So

$$ Vol\biggl(\Pi_{1,@}\big[x',y',H\big]\biggr) =
|x'|\,\mathbb{B}\,\biggl(\frac{Dq\,\lVert \, (y',H) \,   \lVert_{@}}{|x'|}\, \biggr)^2 \stackrel{*}{=}D^2Hq\,\mathbb{B}\,|y|\stackrel{(G.8)}{=}\,\frac{D^2\,\mathbb{B}\,|y|}{L} 
$$ 

where for $\stackrel{*}{=}$ we used $x'=\frac{q\lVert(y',H)\lVert_@^2}{yH}$. \qed

\newpage
\vspace*{0.2cm} \ \\ \textbf{Remark 5.} As we shall see later on, the exact value of $\mathbb{B}$ has no importance for us.  

\ \\ We now introduce an extension for the definitions of cylinders we have been dealing with so far to infinite cylinders.\ \\ For $\vec{u} \in \mathbb{R}^2 $\, and a positive real number R denote 

$$
\boldsymbol{\widehat{\Pi}_{@}\,\biggl(\vec{u}\,,\,R\biggr)} :=  \left\{(\alpha,\vec{y}) \in\mathbb{R}^{3} \mid \,\,\, \lVert \alpha\vec{u}-\vec{y}\, \lVert_{@} \, \leq R \, \right\} $$
and for $(x,y,0)\in G(V)$ define
$$
\boldsymbol{\widehat{\Pi}_{2,@}\bigl[x,y,0\bigr]}\,\,:=
\widehat{\Pi}_@\biggl(\,(\,\dfrac{x\,y}{\Theta_{x}^{2}}\,,\,\dfrac{q\,y}{\Theta_{x}^{2}})\,\,,\,\,\dfrac{Dq\,|y|}{\Theta_{x}} \biggr)
$$ and notice $\widehat{\Pi}_{2,@}\bigl[x,y,0\bigr]$ is an infinite extension of $\Pi_{2,@}\bigl[x,y,0\bigr]$.

\ \\ Now define the set $B_2\subset\mathbb{R}^2$ to be
$$\boldsymbol{B_2}:=\left\{(x,y,0) \in G(V) \,\mid\,\, \bigl(\,\widehat{\Pi}_{2,@}\bigl[x,y,0\,\bigr]\cap G(\Gamma)\subset span_{\mathbb{Z}}(\widetilde{w}_{n-1})\right\}.
$$ 

 \ \\ Informally, $B_2$ is the set of points $(x,y,0)$ such that $G^{-1}\bigl(f_2(x,y,0)\bigr)$ is a "good candidate" for the next best approximation in the sequence $\,\vec{w}_0,\vec{w}_1,...,\vec{w}_{n-1}$ -- without taking into account for the moment the requirement that
 
 \noindent $G^{-1}\bigl(f_2(x,y,0)\bigr)$ needs to be an integer vector. Hence our motivation in what follows is to show that we can find a point $(x,y,0)\in B_2$ such that $G^{-1}\bigl(f_2(x,y,0)\bigr)$ is indeed an integer vector which fulfills properties (2.1)-(2.7) from Theorem 2. In order to do so, firstly we shall investigate the properties of the set $B_2$.

\vspace*{0.3cm}\begin{Lemma}
For all $(x,y,0)\in G(V)$ and all $s\in\{0,1,...,\lfloor D \rfloor\}$ we have that $\,\pm s\, \widetilde{w}_{n-1}\in\widehat{\Pi}_{2,@}\bigl[x,y,0\,\bigr].$
\end{Lemma} 
\ \\ \textit{Proof.}\,\,\,\,\,
$
 \lVert \, \pm sq \,  \bigl(\dfrac{x\,y}{\Theta_{x}^{2}}\,,\,\dfrac{q\,y}{\Theta_{x}^{2}}\bigr) -\vec{0}\,\, \lVert_@\,\, =
\dfrac{s\,q\,|y|}{\Theta_{x}} \leq
\dfrac{Dq\,|y|}{\Theta_{x}}\,,   
$ where the quantity on the right hand side is the radius of the base of $\widehat{\Pi}_{2,@}\bigl[x,y,0\bigr]$. \qed

 \vspace*{0.6cm}

\begin{Lemma}\label{Lemma 4}
$B_2$ is symmetric and star-shaped in $G(V)$ with respect to the origin.\\ I.e. if $(x,y,0)\in B_2\,\cup\partial_{G(V)} B_2$ than for all \, $-1<t<1$ we have that $$ \biggl(tx\,,ty\,,0\biggr)\in B_2. $$
\end{Lemma}

\ \\\textit{Proof.}\,\, By Lemma \ref{Proposition5} property (G.8) we have that 
$$ G(\Gamma)=span_{\mathbb{Z}}\left\{(q,0,0),(q_{n-2}\,,\,L\,,\,0)\right\}=\left\{(sq+iq_{n-2},iL,0)\mid (s,i)\in \mathbb{Z}^2 \right\}.
$$ For the sake of convenience define $\beta=\beta_{s,i}:=sq+iq_{n-2}\,\,and\,\,\, \zeta=\zeta_{s,i}:=iL$ where we omit the dependence on $s,i$ when it is clear from context.\, \\ We have that:
$$
B_2 = \left\{ (x,y,0)\ \middle\vert \begin{array}{l}
   {\forall (s,i)}_{i\neq0}\in\mathbb{Z}^2 \,\,\,\,\,
   \dfrac{Dq\,|y|}{\Theta_{x}} <\,
   \norm{\, \dfrac{\beta_{s,i}\,y\,(x,q)}{\Theta_{x}^2}-(\zeta_{s,i},0)\,
}_@
  \end{array}\right\} $$
\begin{equation}
\label{equation9}
\,\,\,\,\,\,\,\,\,\,\,=
\bigcap_{\substack{(s,i)\in\mathbb{Z}^2   \\
i\neq0}
}
\left\{ (x,y,0)\ \middle\vert \begin{array}{l}
   \dfrac{Dq\,|y|}{\Theta_{x}} <\,
   \norm{\, \dfrac{\beta_{s,i}\,y\,(x,q)}{\Theta_{x}^2}-(\zeta_{s,i},0)\,
   }_@
\end{array}\right\}.\,\,\,\,\,\,\,
\end{equation} So it is enough to show that each set in the intersection above is symmetric and star-shaped with respect to the origin. Assume $(x,y)\in B_2\,\cup\partial\, B_2$. Working on inequality (\ref{equation9}) yields 
 
$$
\dfrac{Dq|y|}{\Theta_{x}}\, <\, \dfrac{1}{\Theta_{x}^2}\,\,
\norm{ \bigl(\beta xy-\zeta\lVert(x,q)\lVert_@^2\,\,,\,\,\beta qy \bigr) }_@
$$

\begin{equation}\label{equation10}
\iff \hspace*{0.5cm}
\dfrac{Dq}{|\beta|}\,\Theta_{x}\,<\, \norm{\biggl(x-\dfrac{\zeta\Theta_{x}^2}{\beta y}\,,\,q\biggr)}_@.
\end{equation} Squaring both sides of the inequality and remembering that by our assumptions\\ $\Theta_{x}=\lVert(x,q)\lVert_@$ and $\norm{(u_1,u_2)}_@=\sqrt{u_1^2+cu_2^2}$ 
we get that equation (\ref{equation10}) holds iff 
\begin{flushleft}
$\hspace*{1cm}
\dfrac{(Dq)^2}{\beta^2}\Theta_{x}^2\,<\,x^2-\dfrac{2\zeta x}{\beta y}(x^2+cq^2)+\dfrac{\zeta^2}{(\beta y)^2}(x^2+cq^2)^2+cq^2$
\end{flushleft}
\begin{flushleft}
$\iff
\dfrac{(Dq)^2}{\beta^2}\,\Theta_{x}^2\,<\,\Theta_{x}^2\biggl(1-\dfrac{2\zeta x}{\beta y}+\dfrac{\zeta^2}{(\beta y)^2}(x^2+cy^2)
\biggr)
$
\end{flushleft}
\begin{flushleft}
$\iff 0\,<\,1-\dfrac{2\zeta x}{\beta y}+\dfrac{\zeta^2}{(\beta y)^2}(x^2+cy^2)-\dfrac{(Dq)^2}{\beta^2}
$
\end{flushleft}
\begin{flushleft}
$\iff
0\,<\,
(\beta y)^2-2\zeta\beta xy +\zeta^2(x^2+cq^2)-(Dqy)^2
$
\end{flushleft}
$$\Longupdownarrow$$
\begin{equation}
\label{equation11}
0\,<\,
(\beta^2-D^2q^2)\,y^2+(\zeta )^2\,x^2+2(-\zeta\beta )\,xy
+c(\zeta q)^2.\,\,\,\,\,\,\,\,\,\,\,\,\,\,\,\,\,\,\,\,\,\,
\end{equation}\ \\
By the fact that an hyperbola is star-shaped with respect to its center, it's enough to notice that equality in (\ref{equation11}) describes an equation of an hyperbola centered at $(0,0)$ and  that the strict inequality part of (\ref{equation11}) corresponds to the connected component of $(0,0)$. Indeed, write
$$A_{yy}:=\beta^2-D^2q^2\,\,\,\,\,\,\,A_{xx}:=\zeta ^2\,\,\,\,\,\,\
A_{xy}:=-\zeta\beta \,\,\,\,\,\,
B_{x}:=0\,\,\,\,\,\,\
B_{y}:=0.\,\,\,\,\,\,\
$$
Equality in (\ref{equation11}) describes an hyperbola if \, $A_{yy}A_{xx}-A_{xy}^2<0$.  (see \cite{Fa})
\begin{flushleft}
$\iff
(\beta^2-D^2q^2)\zeta ^2-(\zeta\beta)^2  \,<\,0
$
\end{flushleft}
\begin{flushleft}
$\iff
-(\zeta Dq)^2  <\,0.$
\end{flushleft} \ \\ The center of the hyperbola $(x_0,y_0)$ is given by
$$x_0=\dfrac{-\bigl(B_xA_{yy}-B_yA_{xy}\bigr)}{A_{yy}A_{xx}-A_{xy}^2} \,\,\,\,\,
y_0=\dfrac{-\bigl(B_yA_{xx}-B_xA_{xy}\bigr)}{A_{yy}A_{xx}-A_{xy}^2}\,\,.
$$

\begin{flushleft}
$\Longrightarrow$\, the center of the hyperbola is at $(0,0)
$.
\end{flushleft}

\ \\ 
If we plug $(x,y)=(0,0)$ in (\ref{equation11}) we get
$0\,<\,c(\zeta q)^2
$
which holds since $0<c$.

\ \\ We deduce each set in the intersection (\ref{equation9}) is the connected component of the plane which contains the origin which is the center of an hyperbola. Summing this up, we deduce that $B_2$ is symmetric and star-shaped with respect to the origin. \\ \qed

\vspace*{0.2cm}
\begin{Lemma}
$B_2$ is an open set inside $G(V)$.
\end{Lemma} \textit{Proof.} \,\,Assume $(x,y,0)\in B_2$, and define
\begin{flushleft}
$F_{x,y}(s,i):=\,
(\beta_{s,i}^2-D^2q^2)\,y^2+\zeta_{s,i} ^2\,x^2+2(-\zeta_{s,i}\beta_{s,i} )\,xy
+c(\zeta_{s,i} q)^2$
\end{flushleft}
\begin{align*}    
=&\,\big((sq+iq_{n-2})^2-D^2q^2\big)\,y^2+(iL )^2\,x^2+2\big(-iL(sq+iq_{n-2}) \big)\,xy
+c(iL q)^2
\\[0.7em]=&\,\big((qy)^2\big)\,s^2+
\big(c(Lq)^2+(Lx)^2-2Lq_{n-2}xy+(q_{n-2}y)^2\big)\,i^2
\\
&+2(qq_{n-2}y^2-Lqxy)si-
(Dqy)^2
\end{align*}

where $\,\beta_{s,i}=\beta,\,\zeta_{s,i}=\zeta$ in the same way as in Lemma \ref{Lemma 4}. 

\ \\ By equation (\ref{equation11}), $(x,y,0)\in B_2$ if and only if $F_{x,y}(s,i)>0$ for all $(s,i)\in\mathbb{Z}^2$ with $i\neq0$.

\ \\ As in the proof of Lemma \ref{Lemma 4}, write \begin{align*}
   A_{ii}:=c(Lq)^2+(Lx)^2-2Lq_{n-2}xy+(q_{n-2}y)^2
   \\A_{ss}:=(qy)^2\,\,\,\,\,\,\,\,\,\,\,\,\,\,\,\,\,\,\,\,\,\,\,\,\,\,\,A_{si}:=qq_{n-2}y^2-Lqxy
\end{align*}
and notice that since $A_{ss}A_{ii}-A_{si}^2=c(Lq^2y)^2
>0$ for any fixed $x,y$, the solutions $s,i$ of $F_{x,y}(s,i)=0$ are on an ellipse. 
\\ As $F_{x,y}(0,0)=-(Dqy)^2<0$, we have that $(x,y,0)\in B_2$ if and only if for all $(s,i)\in\mathbb{Z}^2$ with $i\neq0$ we have that $(s,i)$ lies in the outer part of the ellipse $F_{x,y}(s,i)=0$.

\noindent Continuously varying the values of $x,y$ corresponds to continuously varying the boundary of the ellipse $F_{x,y}(d,i)=0$, which is compact. So for a small enough $\varepsilon$, the condition $F_{x',y'}(s,i)=0$ for some $x',y'$ in the $\varepsilon$-ball around $(x,y)$ can be satisfied only by a finite number of integer pairs. We get that for a small enough $\varepsilon$, we still have that $(s,i)$ lies in the outer part of the ellipse $F_{x+\varepsilon,y+\varepsilon}(s,i)=0$ for all $(s,i)\in\mathbb{Z}^2$ with $i\neq0$. I.e. $(x,y,0)$ is an internal point of $B_2$.\qed

\vspace*{0.2cm}\begin{Lemma}\label{Lemma5}
If $(x,y,0)\in B_2$ then
$(x+k\dfrac{qy}{L},y,0)\in B_2$ 
for all $k\in\mathbb{Z}$. \\
Furthermore, for all $k\in \mathbb{Z}$ the map
$ T(x,y,0)=(x+k\dfrac{qy}{L},y,0)$ acts as the identity map on $span_{\mathbb{R}}(\widetilde{w}_{n-1})$, preserves $G(\Gamma)$ and satisfies 
\begin{equation}
\label{equation12}
T\biggl(\widehat{\Pi}_{2,@}\bigl[x,y,0\,\bigr]\cap G(V)\biggr)=\widehat{\Pi}_{2,@}\bigl[x+k\dfrac{qy}{L},y,0\,\bigr]\cap G(V). \,\,\,\,\,\,\,\,\,\,
\end{equation}

\end{Lemma} \textit{Proof.} \,\,
First notice $T$ preserves the lattice $G(\Gamma)$. Indeed, take $(sq+iq_{n-2},iL,0)\in G(\Gamma)$. \\ So we have \,
$T(sq+iq_{n-2},iL,0)=(sq+iq_{n-2}+k\dfrac{qiL}{L},iL,0)$\\
$ \hspace*{5.08cm}=((s+ik)q+iq_{n-2},iL,0)\in G(\Gamma).$ 

\ \\ Now notice that in order to prove the lemma it's enough to show that equation (\ref{equation12}) holds. Indeed, let $(x,y,0)\in B_2$ and assume by contradiction that $(x+k\dfrac{qy}{L},y,0)\notin B_2$. So by the definition of $B_2$, there exist $(s,i)\in\mathbb{Z}^2$ with $i\neq 0$ such that $(sq+iq_{n-2},iL,0)\in\widehat{\Pi}_{2,@}\bigl[x+k\dfrac{qy}{L},y,0\,\bigr]$.  

\noindent In particular,

$T^{-1}(sq+iq_{n-2},iL,0)\,\in\, T^{-1}\biggl(\widehat{\Pi}_{2,@}\bigl[x+k\dfrac{qy}{L},y,0\,\bigr]\cap G(V)\biggr)$
\ \\ 
$\hspace*{6cm}=\, \widehat{\Pi}_{2,@}\bigl[x,y,0\,\bigr]\cap\, G(V)$
 
 \ \\ where the last equality holds as we assume (\ref{equation12}) holds. \\ As $T$ acts as the identity map on $span_{\mathbb{R}}(\widetilde{w}_{n-1})$ and preserves the lattice $G(\Gamma)$, we deduce that \,$
T^{-1}(sq+iq_{n-2},iL,0) \in
\bigl(\,\widehat{\Pi}_{2,@}\bigl[(x,y,0)\,\bigr]\setminus
\{ s\, \widetilde{w}_{n-1}\,|\,s\in\mathbb{Z}\}\bigr) \cap G(\Gamma)
$\,, contradicting the assumption that $(x,y,0)\in B_2$. 

\ \\  In order to finish we prove equation (\ref{equation12}):  

\ \\ $T\biggl(\widehat{\Pi}_{2,@}\bigl[x,y,0\,\bigr]\,\cap\, G(V)\biggr) = T\biggl( \widehat{\Pi}_@\biggl(\,(\dfrac{x\,y}{\Theta_{x}^{2}}\,,\,\dfrac{q\,y}{\Theta_{x}^{2}})\,\,,\,\,\dfrac{Dq\,\,|y|}{\Theta_{x}} \biggr)  
\cap G(V) \biggr)$ 

\begin{equation*}
=T\biggl(\left\{(\beta,\zeta,0) \in\mathbb{R}^{3} \mid \,\,\, \lVert \beta  \,(\dfrac{x\,y}{\Theta_{x}^{2}}\,,\,\dfrac{q\,y}\Theta_{x}^{2})\,\,
-(\zeta,0)\, \lVert_{@} \, \leq  \, \dfrac{Dq\,\,|y|}{\Theta_{x}} \right\}
\biggr) 
\end{equation*}

$\stackrel{*}{=} T\biggl(\left\{ (\beta,\zeta,0)\in\mathbb{R}^3 \middle\vert 
\begin{array}{l}
   0\,\geq\,
(\beta^2-D^2q^2)\,y^2+\zeta ^2\,x^2+2(-\zeta\beta)\,xy
+c(\zeta q)^2
  \end{array}\right\}
  \biggr)
$

\ \\ $=\left\{ (\beta+\dfrac{kq}{L}\zeta,\zeta,0)\in\mathbb{R}^3 \middle\vert 
\begin{array}{l}
   0\,\geq\,
(\beta^2-D^2q^2)\,y^2+\zeta ^2\,x^2+2(-\zeta\beta)\,xy+c(\zeta q)^2
  \end{array}\right\}  $

\begin{dontbotheriftheequationisoverlong} 
\stackrel{**}{=}\left\{ (\overline{\beta} ,\zeta,0)\in\mathbb{R}^3 \middle\vert 
\begin{array}{l}
   0\,\geq\,
((\overline{\beta}-\dfrac{kq}{L}\zeta)^2-D^2q^2)\,y^2+\zeta ^2\,x^2+2(-\zeta\overline{\beta}+\dfrac{kq}{L}\zeta^2)\,xy
+c(\zeta q)^2
  \end{array}\right\}
  \end{dontbotheriftheequationisoverlong}

\begin{dontbotheriftheequationisoverlong}
=\left\{ (\overline{\beta} ,\zeta,0)\in\mathbb{R}^3 \middle\vert 
\begin{array}{l}
   0\,\geq\,
(\overline{\beta}^2-D^2q^2)\,y^2+\zeta^2(x+\dfrac{kq}{L}\,y)^2+2(-\zeta\overline{\beta})(x+\dfrac{kq}{L}y)\,y+c(\zeta q)^2
  \end{array}\right\}
  \end{dontbotheriftheequationisoverlong}

$= \widehat{\Pi}_{2,@}\bigl[x+k\dfrac{qy}{L},y,0\,\bigr]\cap G(V),$ 

\ \\ where $\stackrel{*}{=}$ follows from the calculations we did in equations (\ref{equation9},\ref{equation11}) (only with changing direction of inequality) and $\stackrel{**}{=}$ holds by plugging\, $\overline{\beta}-\dfrac{kq\zeta}{L}=\beta$. \qed 

$\vspace{0.4cm}$

\begin{Lemma}\label{Lemma6}

Define $\boldsymbol{(x_0,y_0)}:=\sqrt{\frac{4c}{4D^2-1}}(q_{n-2}+\dfrac{1}{2}q,L)$.\ \\ Then we have $(x_0,y_0,0)\in \partial_{G(V)} B_2. $

\end{Lemma} \vspace*{0.2cm}
\textit{Proof.} \,\,We plug $(x,y)=\sqrt{\frac{4c}{4D^2-1}}(q_{n-2}+\dfrac{1}{2}q,L)$ into equation (\ref{equation11}) and get 

\begin{flalign*}
 0\,\leq\,
(\beta^2-D^2q^2)\,(\sqrt{\frac{4c}{4D^2-1}} L)^2+(\zeta )^2\,(\sqrt{\frac{4c}{4D^2-1}} (q_{n-2}+\dfrac{1}{2}q))^2
\\
+2(-\zeta\beta )\,(\sqrt{\frac{4c}{4D^2-1}} (q_{n-2}+\dfrac{1}{2}q))(\sqrt{\frac{4c}{4D^2-1}} L) +c(\zeta q)^2\,. 
\end{flalign*}

\ \\  Remembering $\beta=sq+iq_{n-2}$ and $\zeta=iL$ for $(s,i)\in\mathbb{Z}^2$ with $i\neq0$, we get 

\begin{flalign*}
 0\,\leq\, s^2\bigl(\dfrac{4c}{4D^2-1}\bigr)+i^2\bigl((\dfrac{c}{4D^2-1}+c\bigr)+si\bigl(\dfrac{-4c}{4D^2-1}  \bigr)+1\bigl( \dfrac{4cD^2}{4D^2-1} \bigr)   
\end{flalign*}

$$\Longupdownarrow $$
$$0\,\leq\,s^2+i^2D^2-si-D^2$$

which holds for all $(s,i)\in\mathbb{Z}^2$ with $i\neq0$, with equality for $s=i=1$. \ \\ Hence $\sqrt{\frac{4c}{4D^2-1}}(q_{n-2}+\dfrac{1}{2}q,L,0)\in \partial_{G(V)} B_2.$   

\qed 

\ \\ Throughout previous lemmas we proved properties of the set $B_2$ which are determined only by $\vec{w}_0,...,\vec{w}_{n-2},\vec{w}_{n-1}$. We now continue by stating further properties of $B_2$ which focus on the specific conditions given to us by the $n$th step of the procedure -- the segment $I_n$, the sector $\Omega(\theta,\delta)$, and the congruence property given by $z_n$.

\newpage
\begin{Lemma}\label{Lemma7}
For\,\,$x_0,y_0$ as in Lemma (\ref{Lemma6}), $t_0,\varepsilon\in\mathbb{R}$ and $k\in\mathbb{N}$ define 
 \begin{flalign*}
  \boldsymbol{\bigl(X_{2}[k],Y_2\bigr)}=&\,\,\bigl(X_{2}[k](t_0),Y_2(t_0)\bigr):=t_0\bigl(x_0+\dfrac{kq}{L}y_0\,,\,y_0\bigr)
  \\
\boldsymbol{\Delta_k}=&\,\,\Delta_k(t_0,\varepsilon) :=\dfrac{Y_2\varepsilon}{X_2[k]+\varepsilon}
 \end{flalign*}
and denote \,$\rectangle_{k}=\rectangle_{k}(t_0,\varepsilon,k):=(X_{2}[k]\,,\,Y_2\,,0)-\bigl([\varepsilon,0,0]\times[\Delta_k,0,0]\bigr)$, \\ i.e. $\rectangle_{k}$ is the rectangle in the plane $G(V)$ defined by the following vertices:
\begin{itemize}
    \item $V_1\,:= (X_{2}[k]\,,\,Y_2\,,0)$
    \item $V_2\,:= (X_{2}[k]\,,\,Y_2-\Delta_k\,,\,0)$
    \item $V_3\,:= (X_{2}[k]-\varepsilon\,,\,Y_2-\Delta_k\,,\,0)$
    \item $V_4\,:= (X_{2}[k]-\varepsilon\,,\,Y_2\,,\,0)$
\end{itemize} 
 If $0<t_0<1$ then there exists $0<\varepsilon$ such that for all $k\in\mathbb{N}$ we have\, $\rectangle_{k}\subset B_2$.
\end{Lemma} 

\begin{figure}[H]
    \centering
    \includegraphics[scale=0.7]{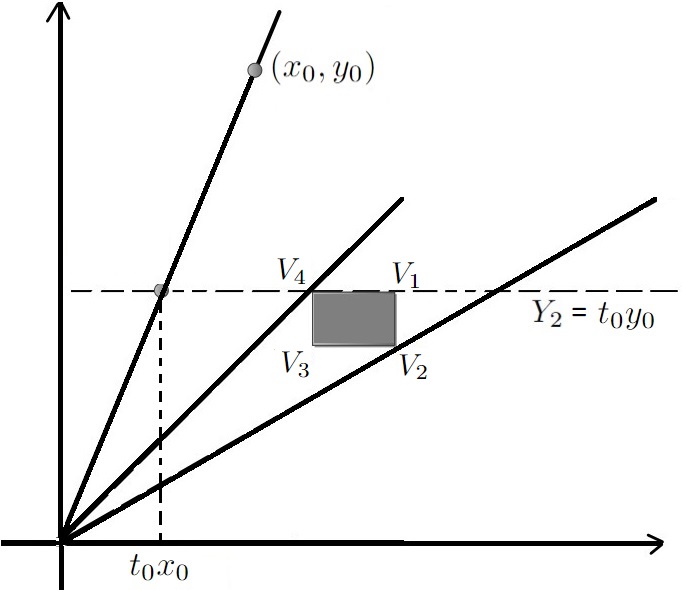}
    \label{fig:my_label}
\end{figure}

\noindent \textit{Proof.}\,\,\, $0<t_0<1$, so by Lemmas \ref{Lemma 4},\ref{Lemma6} we have $t_0\bigl(x_0\,,y_0,\,0\bigr)\in B_2$. 
\\ As $B_2$ is open, there exist $0<\varepsilon$ such that for all $\psi$ with $-\varepsilon<\psi<\varepsilon$ we have that $\bigl(t_0x_0+\psi\,,\,t_0y_0,\,0\bigr)\in B_2$. \\ By Lemma \ref{Lemma5}, for all $k\in\mathbb{N}$ and $-\varepsilon<\psi<\varepsilon$ we have that  $\bigl(t_0x_0+\psi+\dfrac{kq}{L}t_0y_0\,,\,t_0y_0\,,0\bigr)=(X_2[k]+\psi,Y_2,0)\in B_2.$ We use twice again the fact that $B_2$ is star shaped with respect to the origin and get that $\rectangle_{k}\subset B_2$. \qed

\begin{figure}[H]
    \centering
    \includegraphics[scale=0.7]{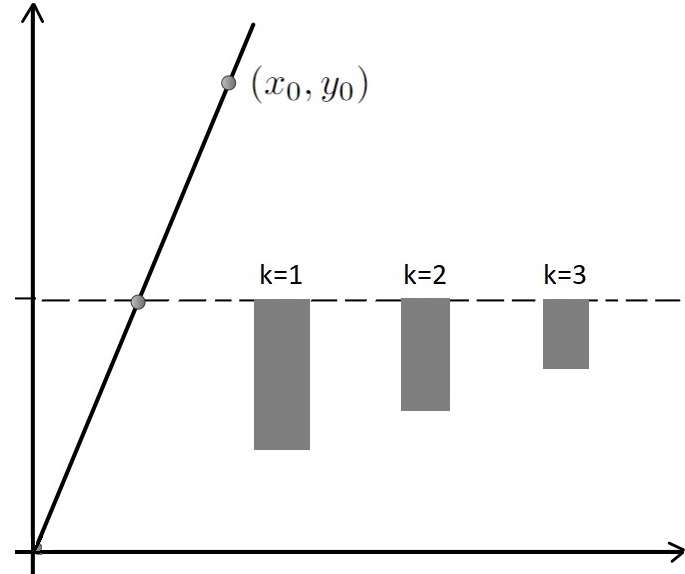}
    \label{fig:my_label}
\end{figure} 

\ \\ Let us observe what happens to $\rectangle_{k}$ under the map $f_2$.

\ \\ As \,$f_2\bigl(x,y,0\bigr)=
(\,\frac{H\,\Theta_{x}^2}{q\,y}\,,\,\frac{H\,x}{q}\,,\,H\,)=(\frac{H(x^2+cq^2)}{qy}\,,\,\frac{H\,x}{q}\,,\,H\,)$ \\ if we fix $y=y_2$ and vary $x$, the line $\{(x,y_2,0)\,|\,x\in\mathbb{R}\}$ maps to the parabola $\{(\frac{H(x^2+cq^2)}{qy_2}\,,\,\frac{H\,x}{q}\,,\,H\,)|\,x\in\mathbb{R}\}$.

\ \\ Similarly, if $x=x_2$ is fixed, the line $\{(x_2,y,0)\,|\,y\in\mathbb{R}\}$ maps to the line \\ $\{(\frac{H(x_2^2+cq^2)}{qy}\,,\,\frac{H\,x_2}{q}\,,\,H\,)\,|\,\,y\in\mathbb{R}\}$.

\ \\ This is demonstrated in the following figure:
\begin{figure}[H]
    \centering
    \includegraphics[scale=0.5]{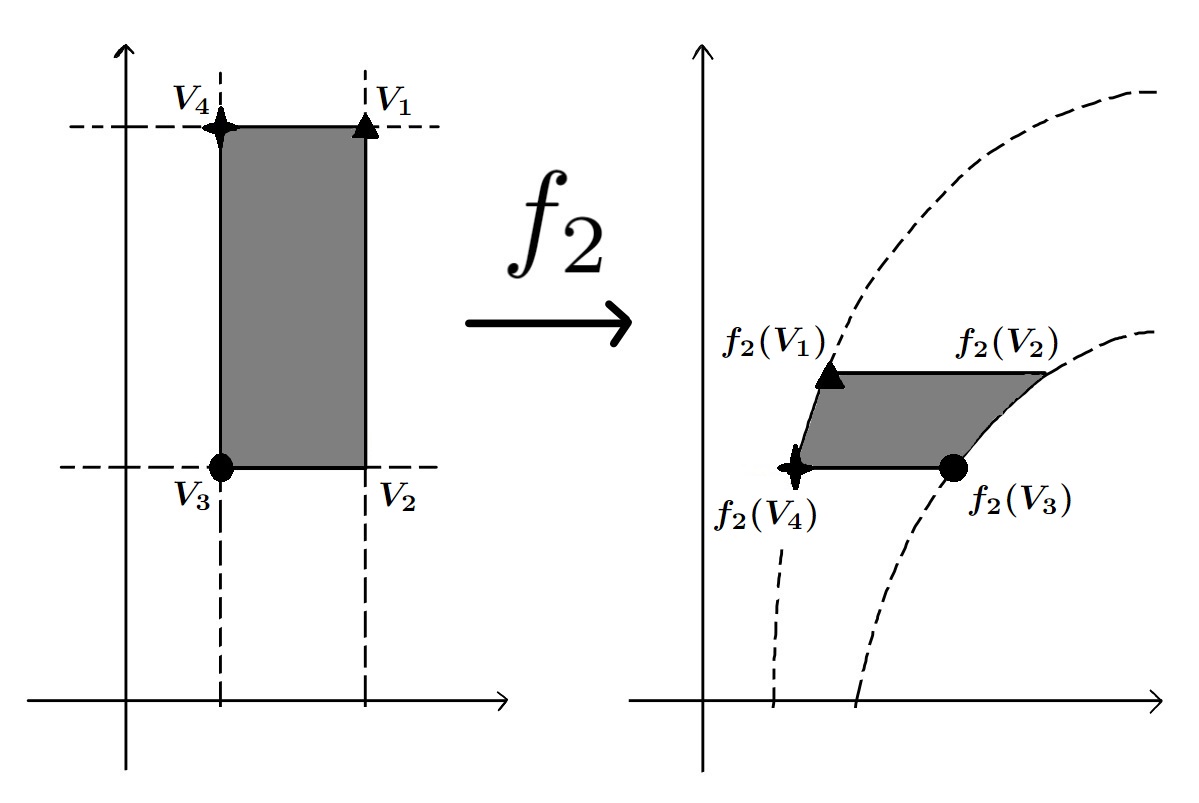}
    \label{fig:my_label}
\end{figure}

\begin{Lemma}\label{Lemma8}
Let $\,0<t_0<1$, $Y_2=t_0y_0$, and let\, $\rectangle_{k}=\rectangle_{k}(t_0,\varepsilon,k)\subset B_2$ be the corresponding rectangles
given to us by Lemma \ref{Lemma7} for some $\varepsilon>0$.

Then the following hold: 
\begin{enumerate}
    \item As $ k\longrightarrow \infty$, the distance between the two branches of the parabolas defining $f_2(\rectangle_{k})$ becomes arbitrarily large.
    \item The distance between the two horizontal lines at the top and bottom of $f_2(\rectangle_{k})$ is fixed and equal $\varepsilon Hq^{-1}$.
    \item The distance between the top horizontal line of $f_2(\rectangle_{k})$ to the top horizontal line of $f_2(\rectangle_{k+1})$ is fixed and equal to $qH^2Y_2$.  
\end{enumerate}
\end{Lemma}
\begin{figure}[H]
    \centering
    \includegraphics[scale=0.5]{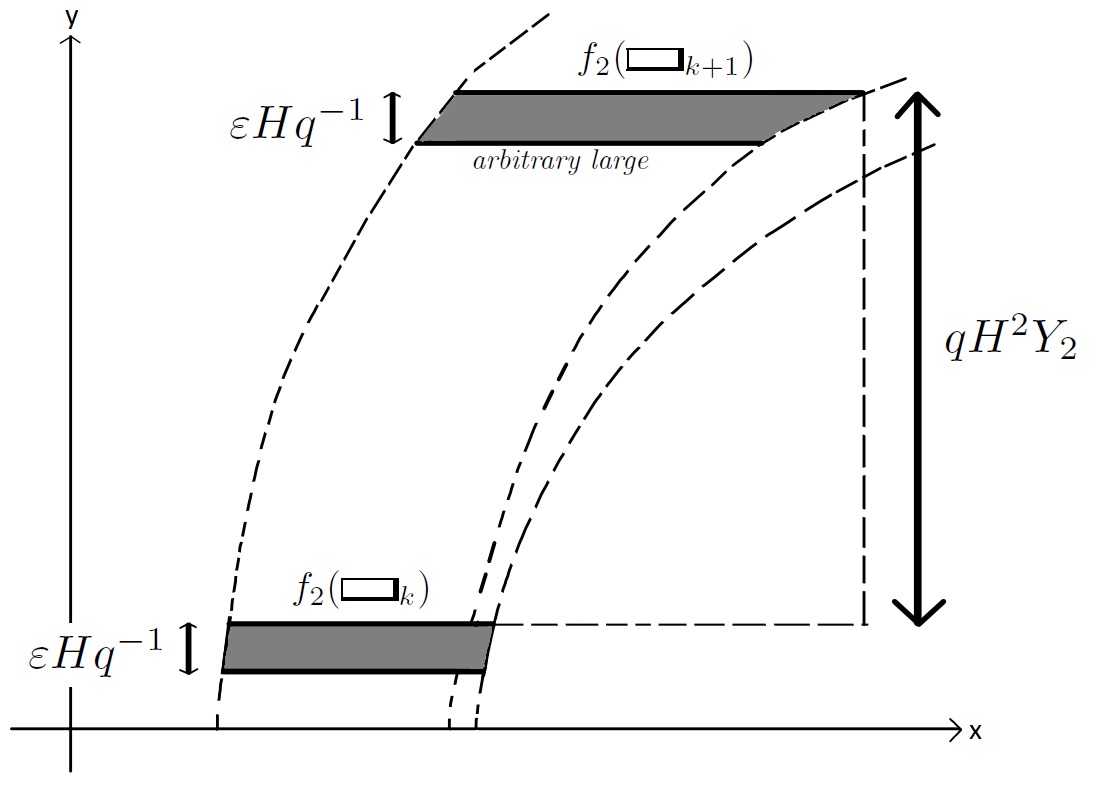}
    \label{fig:my_label}
\end{figure}

\noindent \textit{Proof.}\,\, We prove the claim via direct calculations. For 1, we calculate the difference between the x-value of $f_2(V_3)$ and the x-value of $f_2(V_4)$ and the difference between the x-value $f_2(V_1)$ and the x-value of $f_2(V_2)$ and show both differences diverge to infinity. \ \\ The difference between the x-value of $f_2(V_3)$ and the x-value of $f_2(V_4)$ is equal to \begin{flalign*}
\frac{H\,\lVert(X_2[k]-\varepsilon,q)\lVert_@^2}{q\,(Y_2-\Delta_k)}-\frac{H\,\lVert(X_2[k]-\varepsilon,q)\lVert_@^2}{q\,Y_2}=\frac{\Delta_k \,H\,\lVert(X_2[k]-\varepsilon,q)\lVert_@^2}{q\,(Y_2-\Delta_k)}\\ \sim (X_2[k])^2\,\Delta_k=(X_2[k])^2\,\,\dfrac{Y_2\varepsilon}{X_2[k]+\varepsilon}\xrightarrow {k \to \infty}\infty\,.
\end{flalign*} 

 \newpage \noindent The difference between the x-value of $f_2(V_1)$ and the x-value of $f_2(V_2)$ is equal to 
 \begin{flalign*}
\frac{H\,\lVert(X_2[k],q)\lVert_@^2}{q\,Y_2}-\frac{H\,\lVert(X_2[k],q)\lVert_@^2}{q\,(Y_2-\Delta_k)}=\frac{\Delta_k\, H\,\lVert(X_2[k],q)\lVert_@^2}{q\,(Y_2-\Delta_k)}\\ \hspace*{2cm}\sim (X_2[k])^2\,\Delta_k=(X_2[k])^2\,\,\dfrac{Y_2\varepsilon}{X_2[k]+\varepsilon}\xrightarrow {k \to \infty}\infty\,.
\end{flalign*}

\ \\ For assertion 2, we calculate the difference between the y-value of $f_2(V_1)$ and the y-value of $f_2(V_4)$ and get
\begin{flalign*}
\frac{H}{q}X_2[k]-\frac{H}{q}(X_2[k]-\varepsilon)=\varepsilon Hq^{-1}.
\end{flalign*} 
\ \\ For assertion 3, we calculate the difference between the y-value of $f_2(V_1)$ which corresponds to $f_2(\rectangle_{k+1})$ and the y-value of $f_2(V_1)$ which corresponds to 

\noindent $f_2(\rectangle_{k})$ and get 
\begin{flalign*}
\frac{H}{q}X_2[k+1]-\frac{H}{q}X_2[k]=\frac{H}{q}\bigl(\frac{q}{L}\,Y_2))=H^2q\,Y_2.\qed 
\end{flalign*}

\vspace*{0.2cm}

\begin{Lemma}\label{Lemma9}
Let $\,0<t_0<1$, $Y_2=t_0y_0$, and assume that $H^2qL^{-1}\,Y_2\in\mathbb{R}\setminus\mathbb{Q}$. 
\\ Let\, $\rectangle_{k}=\rectangle_{k}(t_0,\varepsilon,k)\subset B_2$ be the corresponding rectangles
given to us by Lemma \ref{Lemma7} for some $\varepsilon>0$. Let $i_0\in\{0,1,...,m-1\}$ as in (\ref{equation666}). \ \\ Then there exist infinitely many natural numbers $k$  such that $f_2(\rectangle_{k})$ contains at least $m$ elements of $G(\Gamma')$ of the form
\begin{dontbotheriftheequationisoverlong}
\left\{ \bigl(\,(s_k+r)q+(i_0+mi_k)q_{n-2}\,,\,(i_0+mi_k)L\,\,,\,0\,\bigr)+\widetilde{w}_{n-3}\,\,
\middle\vert 
\begin{array}{l}
 r\in\{0,1,...\,,m-1\} \\
\,\,\,\,\,\,\,\,\,\,\,\,s_k,i_k\in\mathbb{Z} 
  \end{array}\right\}
  \end{dontbotheriftheequationisoverlong}
  and with $(s_k+r)q+(i_0+mi_k)q_{n-2} \longrightarrow \infty$ as $ k\rightarrow \infty$.
\end{Lemma}

 \vspace*{0.2cm}
 \noindent \textit{Proof.}\,\,\, By Lemma \ref{Lemma8}, for $k$ large enough the distance between the two branches of the parabolas defining $f_2(\rectangle_{k})$ is greater then $2qm$. Recalling that \\ $G(\Gamma ') = G(\Gamma)+\widetilde{w}_{n-3}= \left\{(sq+iq_{n-2},iL,0)+\widetilde{w}_{n-3}\mid (s,i)\in \mathbb{Z}^2 \right\}
$, in particular we get that if we are able to find $i_k,\,s_k$ such that $(s_kq+(i_0+mi_k)q_{n-2},(i_0+mi_k)L,0)+\widetilde{w}_{n-3}\in f_2(\rectangle_{k})$ for $k$ large enough, then we get $m$ consecutive points of the form $ \bigl((s_k+r)q+(i_0+mi_k)q_{n-2}\,,\,(i_0+mi_k)L\,,\,0\bigr)+\widetilde{w}_{n-3}\,\,
$ with $r\in\{0,1,...,m-1\}$ which are all belong to $f_2(\rectangle_{k})$.

\ \\  Write $\widetilde{w}_{n-3}=(q_{n-3},\widetilde{p}_{n-3})=(q_{n-3},u,H)$. Recalling that the x-values of $f_2(\rectangle_{k})$ lies between $\frac{H(X_2[k]^2+cq^2)}{qY_2}$ and 
$\frac{H(X_2[k]^2+cq^2)}{q(Y_2-\Delta_K)}$, that $\Delta_k\rightarrow 0$, and that $X_2[k]\rightarrow\infty$, we get that in order to prove the lemma it is enough to show we can find infinitely many $i_k\in\mathbb{N}$ such that $(i_0+mi_k)L+u$ \,is between the two horizontal lines at the top and bottom of $f_2(\rectangle_{k})$ for infinitely many $k\in\mathbb{N}$.\ \\ By Lemma \ref{Lemma8} the distance between the top horizontal line of $f_2(\rectangle_{k})$ to the top horizontal line of $f_2(\rectangle_{k+1})$ is fixed and equal to $H^2qY_2$, and the distance between the two horizontal lines at the top and bottom of $f_2(\rectangle_{k})$ is fixed and equal $\varepsilon Hq^{-1}$. I.e. there exist $\alpha,\beta$ such that all top horizontal lines of $f_2(\rectangle_{k})$ are of the form $y=\beta+k(H^2qY_2)$ and all bottom horizontal lines of $f_2(\rectangle_{k})$ are of the form $y=\alpha+k(H^2qY_2)$, and we want to show there exist infinitely many $i_k\in\mathbb{N}$ such that $(i_0+mi_k)L+u\in\big(\,\alpha+k(H^2qY_2)\,,\,\beta+k(H^2qY_2)\,\big)$.
\begin{equation}\label{equation13}
\hspace*{2cm}\iff\,i_k\in\big(\alpha'+k\dfrac{(H^2qY_2)}{mL}\,,\,\beta'+k\dfrac{(H^2qY_2)}{mL}\big) \,\,for
\,some\,\,\alpha',\beta'.
\end{equation} 

\ \\ By the assumptions of the lemma $H^2qY_2L^{-1}\in\mathbb{R}\setminus\mathbb{Q}$, hence $H^2qY_2(mL)^{-1}\in\mathbb{R}\setminus\mathbb{Q}$. So by Kronecker’s density theorem we get infinitely many natural numbers $k,i_k$ such that equation (\ref{equation13}) holds. This finishes the proof of the lemma. \qed

\vspace*{0.3cm}\begin{Lemma}\label{Lemma10}
Let $\widetilde{D}$ such that $\norm{q_{n-2}\vec{v}_{n-1}-\vec{p}_{n-2}}=\widetilde{D}\norm{q_{n-2}\vec{v}_{n-1}-\vec{p}_{n-2}}_@
$. Then there exists $0<t_0<1$ such that the following hold:

\begin{itemize}
    \item $\rectangle_{k}=\rectangle_{k}(t_0,\varepsilon,k)\subset B_2$ for some $\varepsilon>0$ and for all $k\in\mathbb{N}$. 
    \item $H^2qL^{-1}\,Y_2\in\mathbb{R}\setminus\mathbb{Q}$, where $Y_2=t_0y_0$.
    \item $\dfrac{Y_2}{L}\in \widetilde{D}^{-2}\, int\,I_{n-1}^2$.
\end{itemize}
\end{Lemma} 

\ \\ \textit{Proof.}\,\, By Lemma \ref{Lemma7}, there exist $\varepsilon>0$ such that for all $k\in\mathbb{N}$ we have \,$\rectangle_k\subset B_2$. \ \\ We calculate $\dfrac{Y_2}{L}$ and get that $\dfrac{Y_2}{L}=\dfrac{t_0y_0}{L}=\dfrac{t_0\sqrt{\frac{4c}{4D^2-1}}\,L}{L}=t_0\,\sqrt{\frac{4c}{4D^2-1}}$.

\noindent By the conditions of Theorem 2 we have that $I_{n-1}^2\subset[0,M_{\norm{\cdot}}^2]=[0,\sqrt{\frac{4c}{4D^2-1}}]$. 

\noindent So we choose $t_0$ such that $\dfrac{Y_2}{L}\in \widetilde{D}^{-2}\, int\,I_{n-1}^2$\,\,\,and\, $H^2qL^{-1}\,Y_2\in\mathbb{R}\setminus\mathbb{Q}$. \qed 

\vspace*{0.3cm}
\section{Proving Theorem 2 -- Finishing the inductive construction}

Let $0<t_0<1$ be the real number which is given to us by Lemma \ref{Lemma10}, and let $Y_2,X_2[k],\rectangle_k$ be the corresponding points and rectangles in $B_2$ which are given by the lemma for some $\varepsilon>0$. Let $i_0,s_0\in\{0,1,...,m-1\}$ as in (\ref{equation666}), and for convenience write again $G(\vec{w}_{n-3})=(q_{n-3},u,H)$.

\newpage \noindent $H^2q_{n-1}L^{-1}\,Y_2\in\mathbb{R}\setminus\mathbb{Q}$, so by Lemma \ref{Lemma9} there exist infinitely many $k\in\mathbb{N}$ such that $f_2(\rectangle_{k})$ contains at least $m$ elements of $G(\Gamma')$ of the form
\begin{dontbotheriftheequationisoverlong}
\left\{ \bigl(\,(s_k+r)q_{n-1}+(i_0+mi_k)q_{n-2}\,,\,(i_0+mi_k)L\,\,,\,0\,\bigr)+\widetilde{w}_{n-3}\,\,
\middle\vert 
\begin{array}{l}
 r\in\{0,1,...\,,m-1\} \\
\,\,\,\,\,\,\,\,\,\,\,\,s_k,i_k\in\mathbb{Z} 
  \end{array}\right\}
  \end{dontbotheriftheequationisoverlong}
  and with $(s_k+r)q_{n-1}+(i_0+mi_k)q_{n-2} \longrightarrow \infty$ as $ k\rightarrow \infty$. 
  
 \ \\ Our goal is to show that for $k$ large enough we can choose one of these $m$ points to be the vector $\widetilde{w}_n$ we are looking for; i.e. such that \ \\ $\vec{w}_n := G^{-1}(\widetilde{w}_n)=G^{-1}((q_n,\widetilde{p}_n))=(s_k+r)\,\vec{w}_{n-1}+(i_0+mi_k)\,\vec{w}_{n-2}\,+\,\vec{w}_{n-3}$ satisfies properties (2.1)-(2.7).

  \vspace*{0.6cm} \noindent \underline{\textbf{Property \ref{(2.2)}:}}  $G^{-1}$ fixes the x-value, so for $k$ large enough we have that $q_n=(s_k+r)q_{n-1}+(i_0+mi_k)q_{n-2}+q_{n-3}$ is as large as we want.

  \vspace*{0.6cm} \noindent   \underline{\textbf{Property \ref{(2.5)}:}} \,We have that \\ $\lVert \vec{v}_n-\vec{v}_{n-1}\lVert\,\leq D\,\lVert \vec{v}_n-\vec{v}_{n-1}\lVert_@\,\stackrel{(G.10)}{=}D\,\lVert \widetilde{v}_n-\widetilde{v}_{n-1}\lVert_@\,\stackrel{(G.3)}{=}D\,\lVert \widetilde{v}_n\lVert$, 
  \ \\ so in order to prove (2.5) it is enough to show that $\lVert \widetilde{v}_n\lVert\rightarrow0 $
  as $k\rightarrow \infty$. \\ I.e. it is enough to show that $\bigl((i_0+mi_k)L+u\bigr)\bigl((s_k+r)q_{n-1}+(i_0+mi_k)q_{n-2}+q_{n-3}\bigr)^{-1}\rightarrow0$ as $k\rightarrow \infty$. 
  
  \ \\ Recalling that $\bigl(\,(s_k+r)q_{n-1}+(i_0+mi_k)q_{n-2}\,,\,(i_0+mi_k)L\,\,,\,0\,\bigr)+\widetilde{w}_{n-3}\in f_2(\rectangle_k)$, it is enough to show that
  \begin{equation}\label{equation14}
\dfrac{max\,\,\{\,y\,\,\,|\,\,\exists\, t \,\,\,\,such \,\,that \,\,(t,y,H)\in f_2(\rectangle_k)\,\}}{min\,\,\{\,x\,\,|\,\,\exists\, t \,\,\,such \,\,that \,\,(x,t,H)\in f_2(\rectangle_k)\,\}}\longrightarrow\,0\,\,\, as\,\,k \rightarrow\infty. 
  \end{equation}
  \ \\ Indeed, as $Y_2,X_2[k],\Delta_k$ and $\varepsilon$ are positive numbers for all $k\in\mathbb{N}$ we have that the left-hand side of (\ref{equation14}) is equal to

\begin{flalign*}
\dfrac{\dfrac{HX_2[k]}{q_{n-1}}\,}{\dfrac{H\norm{(X_2[k]-\varepsilon,q_{n-1})}_@^2}{q_{n-1}Y_2}}=\dfrac{X_2[k]\,Y_2}{\norm{(X_2[k]-\varepsilon,q_{n-1})}_@^2}\longrightarrow 0\,\,\,as\,\,k\rightarrow\infty.
\end{flalign*}

\vspace*{0.6cm} \noindent
\underline{\textbf{Property \ref{(2.3)}:}} We show that for $1\leq j\leq n$  we have that $\dfrac{R_{n,j}}{R_{n,j-1}}<\dfrac{1}{2}$ \,for $k$ large enough. 

\newpage \noindent For j=n we have 

\begin{flalign*}
\dfrac{R_{n,n}}{R_{n,n-1}}=\dfrac{\norm{q_{n-1}\vec{v}_n-
\vec{p}_{n-1}}}{\norm{q_{n-2}\vec{v}_n-
\vec{p}_{n-2}}}\leq D
\,\dfrac{\norm{q_{n-1}\vec{v}_n-
\vec{p}_{n-1}}_@}{\norm{q_{n-2}\vec{v}_n-
\vec{p}_{n-2}}_@}=
\dfrac{\norm{q_{n-1}\widetilde{v}-
\widetilde{p}_{n-1}
}_@}{\norm{q_{n-2}\widetilde{v}_n-
\widetilde{p}_{n-2}}_@}\ \\ \\=
\dfrac{\norm{q_{n-1}\widetilde{v}_n}_@}{\norm{q_{n-2}\widetilde{v}_n-(L,0)}_@}\longrightarrow\dfrac{0}{\norm{(L,0)}}=0\,\,\,as\,k\rightarrow\infty.
\end{flalign*}

\ \\If $1\leq j\leq n-1$ then since  $\vec{v}_n\longrightarrow\vec{v}_{n-1}$ (property (2.5)) we have

\begin{flalign*}
\dfrac{R_{n,j}}{R_{n,j-1}}=\dfrac{\norm{q_{j-1}\vec{v}_n-
\vec{p}_{j-1}}}{\norm{q_{j-2}\vec{v}_n-
\vec{p}_{j-2}}}\longrightarrow \dfrac{\norm{q_{j-1}\vec{v}_{n-1}-
\vec{p}_{j-1}}}{\norm{q_{j-2}\vec{v}_{n-1}-
\vec{p}_{j-2}}}=\dfrac{R_{n-1,j}}{R_{n-1,j-1}}<\dfrac{1}{2}
\end{flalign*}
\ \\where the last inequality holds by our inductive assumption.

\vspace*{0.6cm} \noindent
\underline{\textbf{Property \ref{(2.6)}:}} We show that for $1\leq j\leq n$  we have that $\,\, \vec{\psi}_{n,j}\in \, int\, \Omega(\theta,\delta)$ for $k$ large enough.  

\ \\ For $j=n$ we have: 

\begin{flushleft}
$\vec{\psi}_{n,n}=q_{n-1}\vec{v}_{n}-\vec{p}_{n-1}=\dfrac{1}{q_n}\,(q_{n-1}\vec{p}_n-q_n\vec{p}_{n-1})$
\end{flushleft}
\begin{flushleft}
$=\dfrac{1}{q_n}\,\biggl(q_{n-1}\big((s_k+r)\vec{p}_{n-1}+(i_0+mi_k)\vec{p}_{n-2}+\vec{p}_{n-3}\big)$
\\
$\hspace*{3.4cm}-\big((s_k+r)q_{n-1}+ (i_0+mi_k)q_{n-2}+q_{n-3}\big)\,\vec{p}_{n-1}\biggr)$
\end{flushleft}
\begin{equation}\label{equation15}
=\dfrac{1}{q_n}\,\biggl(
(i_0+mi_k)\,(q_{n-1}\vec{p}_{n-2}-q_{n-2}\vec{p}_{n-1})+q_{n-1}\vec{p}_{n-3}-q_{n-3}\vec{p}_{n-1} \biggr)  
\end{equation}\ \\ As $k\rightarrow\infty$ we have that $i_k\rightarrow\infty$, and so the direction of $\vec{\psi}_{n,n}$ converges to the direction of  $q_{n-1}\vec{p}_{n-2}-q_{n-2}\vec{p}_{n-1}$
$=-q_{n-1}(\vec{\psi}_{n-1,n-1})$, i.e. to the direction of $-\vec{\psi}_{n-1,n-1}$.

\ \\ By our inductive assumption $\vec{\psi}_{n-1,n-1}\in int\, \Omega(\theta,\delta)$. So for $k$ large enough we have that $ \vec{\psi}_{n,n}\in \, int\, \Omega(\theta,\delta)$ and that $\vec{\psi}_{n,n}$ and $\vec{\psi}_{n-1,n-1}$ lie on opposite quadrants.

\ \\ If $1\leq j\leq n-1$ then as $\vec{v}_n\longrightarrow \vec{v}_{n-1}$ we have that $\vec{\psi}_{n,j}\longrightarrow \vec{\psi}_{n-1,j}$. \ \\ By our inductive assumption $\vec{\psi}_{n-1,j}\in int\, \Omega(\theta,\delta)$. Furthermore, for $2\leq j\leq{n-1}$ we have that $\vec{\psi}_{n-1,j}$ and $\vec{\psi}_{n-1,j-1}$ lie on opposite quadrants. So for $k$ large enough we get property (2.6).

\newpage \noindent
\underline{\textbf{Property \ref{(2.4)}:}} We show that for $1\leq j\leq n$  we have that \\ $ q_{j}\,\,   \lVert q_{j-1}\vec{v}_n-\vec{p}_{j-1} \lVert^{2}\,\, \in \, int\,I_{j-1}^2$ for $k$ large enough.

\ \\ We start with the case $j=n$.
\ \\ By equation (\ref{equation8}) and Lemma \ref{Lemma3}, if $(x,y,0)\in B_2$ and $(q_n,\vec{p}_n)=G^{-1}
(q_n,\widetilde{p}_n)=G^{-1}
\big(f_2(x,y,0)\big) 
$ then we have that
$$q_n\norm{q_{n-1}\vec{v}_n-\vec{p}_{n-1}}^2=\dfrac{\mid y\mid D^2\,\mathbb{B}\,L^{-1}}{\mathbb{B}\,(D\Lambda)^2\,}=\dfrac{\mid y\mid }{\Lambda^2\,L}
$$ where $\Lambda$ is a constant satisfying 
$\norm{q_{n-1}\vec{v}_n-\vec{p}_{n-1}}_@=\Lambda\norm{q_{n-1}\vec{v}_n-\vec{p}_{n-1}}$ with $\frac{1}{D}\leq \Lambda\leq1$. 

\ \\ By property (2.6), $\Lambda\rightarrow \widetilde{D}^{-1}$ as $k\rightarrow\infty$. \\ For every $k$, the top horizontal line defining $\rectangle_k$ is equal to $y=Y_2$, and the bottom horizontal line defining $\rectangle_k$ is equal to $y=Y_2-\Delta_k$ with $\Delta_k\rightarrow 0$. \ \\ So if $(q_n,\widetilde{p}_n)\in f_2(\rectangle_{k})$ then  we have that $q_n\norm{q_{n-1}\vec{v}_n-\vec{p}_{n-1}}^2\rightarrow\dfrac{\widetilde
{D}^2Y_2}{L}$ \\ as $k\rightarrow\infty$. 

\ \\ As we chose $Y_2$ to satisfy $\dfrac{Y_2}{L}\in \widetilde{D}^{-2}\,int\, I_{n-1}^2$ (Lemma \ref{Lemma10}), we get in total that for $k$ large enough\, $q_n\norm{q_{n-1}\vec{v}_n-\vec{p}_{n-1}}^2\in int\,I_{n-1}^2$.

\ \\ If $1\leq j\leq n-1$ then as $\vec{v}_n\rightarrow \vec{v}_{n-1}$ we have that $ q_{j}\,\,   \lVert q_{j-1}\vec{v}_n-\vec{p}_{j-1} \lVert^{2}\longrightarrow q_{j}\,\,   \lVert q_{j-1}\vec{v}_{n-1}-\vec{p}_{j-1} \lVert^{2}$. \\By our inductive assumption $q_j\,\,\lVert q_{j-1}\vec{v}_{n-1}-\vec{p}_{j-1} \lVert^{2}\in\,int\,I_{j-1}^2$, and so for $k$ large enough we get property (2.4).

\vspace*{0.6cm} \noindent
\underline{\textbf{Property \ref{(2.1)}:}}\,  We start with the following observation:

\begin{Lemma}\label{Lemma14}
Let $(\vec{a}_k)_{k=0}^{\infty
}$ and $(\vec{b}_k)_{k=0}^{\infty
}$ be two sequences in $\mathbb{R}^d$, let $\norm{
\cdot
}$ be an arbitrary norm, and assume that:
\begin{enumerate}
    \item $\lim\limits_{k\rightarrow\infty}\vec{a}_k=\vec{t}$ \,for $\vec{t}\neq\vec{0}$.
    \item \,$\vec{b}_k\neq\vec{0}$ for all $k$.
    \item The direction of $\vec{b}_k$ converges to the direction of $\vec{t}$.
\end{enumerate}

Then for $k$ large enough we have that $\lVert \vec{a}_k \lVert\, < \lVert \vec{a}_k+\vec{b}_k \lVert.$

\end{Lemma} 

\vspace*{0.2cm} \noindent 
\textit{Proof.} \, Let $(s_k)_{k=0}^{\infty}$ be the sequence of positive numbers defined by $s_k:=\lVert\vec{b}_k\lVert\cdot \lVert\vec{a}_k\lVert^{-1}$, and let $\vec{t}_k$ be the sequence of vectors defined by
$\vec{t}_k:=\vec{b}_k-s_k\vec{a}_k$. So particularly we have that \, $\dfrac{1}{s_k}\vec{t}_k=\dfrac{1}{s_k}\vec{b}_k-\vec{a}_k$. \ \\ By assumption 3, the direction of $\dfrac{1}{s_k}\vec{b}_k$ and the direction $-\vec{a}_k$ converge oppositely, both having the same length (w.r.t $\norm{\cdot}$), hence $\lim\limits_{k\rightarrow\infty}\dfrac{1}{s_k}\vec{t}_k=\vec{0}$. \\ So we have
\begin{flalign*}
\lVert \vec{a}_k+\vec{b}_k \lVert \,=\, \lVert \vec{a}_k + s_k\vec{a}_k+\vec{t}_k  \lVert \,\geq\, (1+s_k)\lVert \vec{a}_k \lVert-\lVert
\vec{t}_k\lVert\,>\,\lVert\vec{a}_k\lVert,
\end{flalign*}
where the last inequality holds for $k$ large enough iff \,$\dfrac{1}{s_k}\lVert\vec{t}_k\lVert\,<\,\lVert\vec{a}_k\lVert$, which holds for $k$ large enough as $\lim\limits_{k\rightarrow\infty}\dfrac{1}{s_k}\vec{t}_k=\vec{0}$ and $\lim\limits_{k\rightarrow\infty}\vec{a}_k=\vec{t}$ \,for $\vec{t}\neq\vec{0}$.  \\ \qed

\ \\ Now we show that if $k$ is large enough then the followings hold:

\begin{itemize}
    \item $\Pi_{n,n} \cap \mathbb{Z}^3 = 
    \{\vec{w}_n, \vec{w}_{n-1}, \vec{w}_n-\vec{w}_{n-1},\vec{0} \}$ \,with\,  $int\,\Pi_{n,n} \cap \mathbb{Z}^3 =\varnothing$ 
    
    \item $\Pi_{n,j} \cap \mathbb{Z}^3 = 
    \{\vec{w}_j, \vec{w}_{j-1},\vec{0}\} $\,\,with\,\,\,$  int\,\Pi_{n,j} \cap \mathbb{Z}^3 =\varnothing$. 
\end{itemize}

 \vspace*{0.2cm}
 We start with the case $j=n$. \ \\ As $\Pi_{n,n}=\left\{(\alpha,\vec{y}) \in\mathbb{R}^{3} \mid 0 \leq \alpha \leq q_n \,\,\, and \,\,\, \lVert \alpha\vec{v}_n-\vec{y}\, \lVert \, \leq \norm{q_{n-1}\vec{v}_n-\vec{p}_{n-1}} \, \right\} $, \\we immediately get that $\{ \vec{w}_n,\vec{w}_{n}-\vec{w}_{n-1},\vec{w}_{n-1},\vec{0}\,
\}\subset\partial\Pi_{n,n}$. 
\ \\Furthermore, for all $s\geq2$ we have $\lVert sq_{n-1}\vec{v}_n-s\vec{p}_{n-1}\, \lVert \, > \norm{q_{n-1}\vec{v}_n-\vec{p}_{n-1}}$, hence for all $s\geq2$ we have that $s\vec{w}_{n-1}\notin\Pi_{n,n}$.

\ \\ Now notice that the map $s\vec{w}_{n-1}+i\vec{w}_{n-2}\longmapsto\vec{w}_n-(s\vec{w}_{n-1}+i\vec{w}_{n-2})$ is a bijection from $\Gamma$ to $\Gamma'$ such that 
\ \\
$\norm{(sq_{n-1}+iq_{n-2})\vec{v}_n-(s\vec{p}_{n-1}+\vec{p}_{n-2)}}=$
\\$\hspace*{3.5cm}\norm{(q_n-(sq_{n-1}+iq_{n-2}))\vec{v}_n-(\vec{p}_n-(s\vec{p}_{n-1}+\vec{p}_{n-2}))}$. 

\ \\ So we also have that $\vec{w}_n-s\vec{w}_{n-1}\notin\Pi_{n,n}$ for all $s\geq2$.  

\ \\ Now notice that by definitions we have $\Pi_{n,n}\subseteq\Pi_@(\vec{v}_{n},q_{n},DR_{n,n})$. Recalling property (G.\ref{G.11}) we get that it is now enough to show that \,\,$$int\,\Pi_@(\widetilde{v}_{n},q_{n},DR_{n,n}) \cap G(\mathbb{Z}^3)\setminus
\{ s\, \widetilde{w}_{n-1}\,,\widetilde{w}_n-s\widetilde{w}_{n-1}\,|\,s\in\mathbb{N}\}=\varnothing.$$ 

\ \\ Recalling that by equation (\ref{equation7}) we have $\Pi_@(\widetilde{v}_{n},q_{n},DR_{n,n})=\Pi_{1,@}\big[q_n,\widetilde{p}_n\big]$, we get that it is enough to show that 

\,\,$$int\,\Pi_{1,@}\big[q_n,\widetilde{p}_n\big] \cap G(\mathbb{Z}^3)\setminus
\{ s\, \widetilde{w}_{n-1}\,,\widetilde{w}_n-s\widetilde{w}_{n-1}\,|\,s\in\mathbb{N}\}=\varnothing.$$

\ \\ For all $k$ we have that $\rectangle_k\subset B_2$, so by the definition of $B_2$ and Lemma \ref{Lemma2} we have that $int\,\Pi_{1,@}\big[q_n,\widetilde{p}_n\big]\cap G(\Gamma)\setminus
\{ s\, \widetilde{w}_{n-1}\,|\,s\in\mathbb{N}\}=\varnothing$. 

\ \\ Now notice that the map $s\widetilde{w}_{n-1}+i\widetilde{w}_{n-2}\longmapsto\widetilde{w}_n-(s\widetilde{w}_{n-1}+i\widetilde{w}_{n-2})$ is a bijection from $G(\Gamma)$ to $G(\Gamma')$ such that

$\norm{(sq_{n-1}+iq_{n-2})\widetilde{v}_n-(s\widetilde{p}_{n-1}+\widetilde{p}_{n-2)}}=$

$\hspace*{3.5cm}\norm{(q_n-(sq_{n-1}+iq_{n-2}))\widetilde{v}_n-(\widetilde{p}_n-(s\vec{p}_{n-1}+\widetilde{p}_{n-2}))}$. 

\ \\ So we also have that $int\,\Pi_{1,@}\big[q_n,\widetilde{p}_n\big]\cap G(\Gamma')\setminus
\{ \widetilde{w}_n-s\, \widetilde{w}_{n-1}\,|\,s\in\mathbb{N}\}=\varnothing$.

\ \\ Since $\{\vec{w}_{n-1},\vec{w}_{n-2},\vec{w}_{n-3}\}$ is a basis of $\mathbb{Z}^3$, $V=span_{\mathbb{R}}\,\{\vec{w}_{n-1},\vec{w}_{n-2}\}$ and \\$V'=V+\vec{w}_{n-3}$, we know that $G(\mathbb{Z}^3)\cap\{(x,y,t)\,|\,0<t<H\,\}=\varnothing$.

\ \\ $\Pi_{1,@}\big[q_n,\widetilde{p}_n\big]$ is an elliptical cylinder (w.r.t $\norm{\cdot}_@$) such that the axis of the cylinder is the line segment which connects the origin and $(q_n,\widetilde{p}_n)$. The radius of the cylinder is equal to $D\,\norm{q_{n-1}\widetilde{v}_n-\widetilde{p}_{n-1}}=D\,\norm{q_{n-1}\widetilde{v}_n}$, which converges to $0$ as $k\rightarrow\infty$. So for $k$ large enough, $int\, \Pi_{1,@}\big[q_n,\widetilde{p}_n\big] \cap G(\mathbb{Z}^3)\subset G(\Gamma)\cup G(\Gamma')$, which finishes the case $j=n$.

\ \\ For the case $j=n-1$, as \\ $\Pi_{n,n-1}=\left\{(\alpha,\vec{y}) \in\mathbb{R}^{3} \mid 0 \leq \alpha \leq q_{n-1} \,\,\, and \,\,\, \lVert \alpha\vec{v}_n-\vec{y}\, \lVert \, \leq \norm{q_{n-2}\vec{v}_n-\vec{p}_{n-2}} \, \right\} $, \\we immediately get that $\vec{w}_{n-2}
\in\partial\Pi_{n,n-1}$. By (2.3), for $k$ large enough we have $\norm{q_{n-1}\vec{v}_n-\vec{p}_{n-1}}<\norm{q_{n-2}\vec{v}_n-\vec{p}_{n-2}}$, and so for $k$ large enough we also have $\vec{w}_{n-1}
\in\partial\Pi_{n,n-1}$. 

\ \\ By the inductive assumption $\Pi_{n-1,n-1} \cap \mathbb{Z}^3 = 
    \{\vec{w}_{n-1}, \vec{w}_{n-2},\vec{w}_{n-1}- \vec{w}_{n-2},\vec{0}\} $. 
\ \\ Since $\vec{v}_n\rightarrow\vec{v}_{n-1}$ as $k\rightarrow\infty$,
we have that $\Pi_{n,n-1}\rightarrow\Pi_{n-1,n-1}$ as $k\rightarrow\infty$.
\ \\ So in order to finish, we only need to show that for $k$ large enough \\ $\vec{w}_{n-1}-\vec{w}_{n-2}\notin \,\Pi_{n,n-1}$.

\begin{flalign*}
\iff 
\norm{q_{n-2}\vec
{v}_n-\vec{p}_{n-2}}
&<
\norm{(q_{n-1}-q_{n-2})\vec{v}_n-(\vec{p}_{n-1}-\vec{p}_{n-2})}
\\&=\norm{q_{n-1}\vec{v}_n-\vec{p}_{n-1}-(q_{n-2}\vec{v}_n-\vec{p}_{n-2})}
\\&=\norm{q_{n-2}\vec{v}_n-\vec{p}_{n-2}-(q_{n-1}\vec{v}_n-\vec{p}_{n-1})}.
 \end{flalign*}
 
\vspace*{0.2cm}
\noindent That is, we need to show that for $k$ large enough
\begin{equation}\label{equation16}
\norm{\vec{\psi}_{n,n-1}}<\norm{\vec{\psi}_{n,n-1}-\vec{\psi}_{n,n}}.  
\end{equation}
\ \\ By equation (\ref{equation15}), the direction of $\psi_{n,n}$ converges to the direction of $-\psi_{n-1,n-1}$ as $k\rightarrow\infty$. For the direction of $\psi_{n,n-1}$, we have that
\begin{flushleft}
$\vec{\psi}_{n,n-1}=q_{n-2}\vec{v}_{n}-\vec{p}_{n-2}=\dfrac{1}{q_n}\,(q_{n-2}\vec{p}_n-q_n\vec{p}_{n-2})$
\end{flushleft}
\begin{flushleft}
$=\dfrac{1}{q_n}\,\biggl(q_{n-2}\big((s_k+r)\vec{p}_{n-1}+(i_0+mi_k)\vec{p}_{n-2}+\vec{p}_{n-3}\big)$ \\
$\hspace*{4.6cm}-\big((s_k+r)q_{n-1}+(i_0+mi_k)q_{n-2}+q_{n-3}\big)\,\vec{p}_{n-2}\biggr)$
\end{flushleft}
\begin{flushright}
$=\dfrac{1}{q_n}\,\biggl(
(s_k+r)\,(q_{n-2}\vec{p}_{n-1}-q_{n-1}\vec{p}_{n-2})+q_{n-2}\vec{p}_{n-3}-q_{n-3}\vec{p}_{n-2} \biggr)$. 
\end{flushright}
\ \\ As $k\rightarrow\infty$ we have that $s_k\rightarrow\infty$, and so the direction of $\vec{\psi}_{n,n-1}$ converges to the direction of $q_{n-2}\vec{p}_{n-1}-q_{n-1}\vec{p}_{n-2}=q_{n-1}(\vec{\psi}_{n-1,n-1})$, i.e. to the direction of $\vec{\psi}_{n-1,n-1}$.

\ \\ So we have that $\vec{\psi}_{n,n-1}=q_{n-2}\vec
{v}_n-\vec{p}_{n-2}$ converges to $
\vec{\psi}_{n-1,n-1}\neq
\vec{0}$ as $k\rightarrow\infty$, 
the direction of $\vec{\psi}_{n,n-1}$ and the direction $-\vec{\psi}_{n,n}$ are both converging to the same direction as $k\rightarrow\infty$, both being sequences of non-zero vectors. We deduce by Lemma \ref{Lemma14} that for $k$ large enough equation (\ref{equation16}) holds.

\ \\ Lastly, we deal with the case $1\leq j< n-1$. \ \\ As $\Pi_{n,j}=\left\{(\alpha,\vec{y}) \in\mathbb{R}^{3} \mid 0 \leq \alpha \leq q_j \,\,\, and \,\,\, \lVert \alpha\vec{v}_n-\vec{y}\, \lVert \, \leq \norm{q_{j-1}\vec{v}_n-\vec{p}_{j-1}} \, \right\} $, \\we immediately get that $ \vec{w}_{j-1}
\in\partial\Pi_{n,j}$. By (2.3), for $k$ large enough we have $\norm{q_{j-1}\vec{v}_n-\vec{p}_{j-1}}<\norm{q_{j-2}\vec{v}_n-\vec{p}_{j-2}}$, and so for $k$ large enough we also have $\vec{w}_{j}
\in\partial\Pi_{n,j}$. 

\ \\ By the inductive assumption $\Pi_{n-1,j} \cap \mathbb{Z}^3 = 
\{\vec{w}_j, \vec{w}_{j-1},\vec{0}\}$. Once again, since $\vec{v}_n\rightarrow\vec{v}_{n-1}$ as $k\rightarrow\infty$,
we have that $\Pi_{n,j}\rightarrow\Pi_{n-1,j}$ as $k\rightarrow\infty$. So by the inductive assumption, for k large enough we have $\Pi_{n,j} \cap \mathbb{Z}^3 =\{\vec{w}_j,\vec{w}_{j-1},\vec{0}\}$ and $int\,\Pi_{n,j} \cap \mathbb{Z}^3 =\varnothing$, and so in total we get property (2.1).

\vspace*{0.6cm}\noindent
\underline{\textbf{Property \ref{(2.7)}:}}
 \, Throughout the current chapter we showed that we have infinitely many points of the form \begin{equation}\label{equation18}
  \vec{w}_n=(s_k+r)\,\vec{w}_{n-1}+(i_0+mi_k)\,\vec{w}_{n-2}\,+\,\vec{w}_{n-3}   
 \end{equation} with $r\in\{0,1,...\,,m-1\}$, and all of these $m$ points satisfy properties (2.1)-(2.6) assuming $k$ is large enough. 

 \ \\ In particular, $q_n=(s_k+r)\,q_{n-1}+(i_0+mi_k)\,q_{n-2}\,+\,q_{n-3}$. 
  
\ \\ Recall that by (\ref{equation666}), our goal is to have $q_n\equiv s_0z_{n-1}+i_0z_{n-2}+z_{n-3}\mod{m}$. 

\ \\ \noindent By the inductive assumption we have
$\vspace*{0.2cm}$
$\Bigg\{
\begin{array}{l}
q_{n-1}\equiv z_{n-1}\mod m
   \vspace*{0.1cm}\\
q_{n-2}\equiv z_{n-2}\mod m
 \vspace*{0.1cm}\\
q_{n-3}\equiv z_{n-3}\mod m.
\end{array}$

\ \\ So for each $k$, we choose $r=r_k$ such that $s_k+r\equiv s_0\mod m$ \,and get infinitely many points which satisfy property (2.7).

\ \\ As $\vec{w}_n$ satisfies (\ref{equation18}) and $\{\vec{w}_{n-1},\vec{w}_{n-2},\vec{w}_{n-3}\}$ a basis of $\mathbb{Z}^3$, we also have that $\{\vec{w}_n,\vec{w}_{n-1},\vec{w}_{n-2}\}$ is a basis of $\mathbb{Z}^3$. 
  
\ \\ Furthermore, as $\vec{w}_n$ is the $(n+1)$th best approximation vector of $\vec{v}=\lim\limits_{n\rightarrow\infty} \vec{v}_n$, and as we have infinitely many candidates for $\vec{w}_n$ which satisfy properties (2.1)-(2.7) at each step, we get by uniqueness of the best approximation vectors sequence (starting from some index) continuum many vectors 
  which satisfy Theorem 2. This finishes the proof. \ \\  \qed

\newpage
\bibliographystyle{alpha}
\bibliography{ref}

\end{document}